\documentclass[final,hidelinks,onefignum,onetabnum,reqno,dvipsnames]{siamart250211}

\usepackage{amsfonts} \usepackage{graphicx} \usepackage{booktabs}
\usepackage{subcaption}
\usepackage{epstopdf}
\usepackage{algorithm}
\usepackage{algpseudocodex}
\usepackage{comment}
\usepackage{multirow}
\usepackage{longtable}
\usepackage{xcolor}
\usepackage{xparse}
\ifpdf
	\DeclareGraphicsExtensions{.eps,.pdf,.png,.jpg}
\else
	\DeclareGraphicsExtensions{.eps}
\fi

\usepackage{cite}
\usepackage[shortlabels]{enumitem}

\makeatletter
\patchcmd{\ALG@step}{\addtocounter{ALG@line}{1}}{\refstepcounter{ALG@line}}{}{}
\newcommand{\ALG@lineautorefname}{Step}
\makeatother

\newsiamremark{remark}{Remark}
\newsiamremark{hypothesis}{Hypothesis}
\crefname{hypothesis}{Hypothesis}{Hypotheses}
\newsiamthm{claim}{Claim}

\headers{ROED of infinite-dimensional Bayesian nonlinear inverse problems}{A. Chowdhary, A. Attia, and A. Alexanderian}

\title{
	Robust optimal experimental design of infinite-dimensional Bayesian nonlinear inverse problems%
	\thanks{
		Submitted to the editors \today.
		\funding{
			The work of A.~Chowdhary and A.~Alexanderian was supported in part by US National
			Science Foundation grants DMS \#2111044.  The work of A.~Alexanderian was also
			supported in part by US National Science Foundation grants DMS \#1745654.
			The work of Ahmed Attia is supported by the U.S. Department of Energy, Office of
			Science, Office of Advanced Scientific Computing Research (ASCR), ASCR Applied
			Mathematics Base Program, and Scientific Discovery through Advanced Computing
			(SciDAC) Program through the FASTMath Institute under contract number
			DE-AC02-06CH11357 at Argonne National Laboratory.
		}
	}
}

\author{
	Abhijit Chowdhary\thanks{
		Department of Mathematics, North Carolina State University, Raleigh, NC
		(\email{achowdh2@ncsu.edu}, \email{alexanderian@ncsu.edu})
	}
	\and Ahmed Attia\thanks{
		Mathematics and Computer Science Division, Argonne National Laboratory, Lemont, IL (\email{attia@mcs.anl.gov})
	}
	\and Alen Alexanderian\footnotemark[2]
}

\ifpdf
	\hypersetup{
		pdftitle={ Robust optimal design of large-scale Bayesian nonlinear inverse problems},
		pdfauthor={A. Chowdhary, A. Attia, and A. Alexanderian}
	}
\fi

\usepackage{amsmath, amssymb}
\usepackage{commath} 
\usepackage{bm}
\usepackage{dsfont}
\usepackage{mathtools}
\usepackage{mathrsfs}
\usepackage{physics}
\usepackage{stmaryrd}
\usepackage{amsopn}
\usepackage[normalem]{ulem}
\usepackage{adjustbox}

\allowdisplaybreaks

\usepackage{tcolorbox}
\tcbuselibrary{skins} 
\tcbuselibrary{theorems} 
\tcbuselibrary{breakable} 
\usetikzlibrary{shadows.blur} 
\newtcbtheorem[crefname={problem}{problems}]{problem}{Problem}{%
	enhanced jigsaw,
	colback=gray!20!white,
	colframe=gray!80!black,
	size=small,
	boxrule=1pt,
	title=\textbf{Problem},
	halign title=flush center,
	coltitle=black,
	breakable,
	drop shadow=black!50!white,
	attach boxed title to top left={xshift=0.5cm,yshift=-\tcboxedtitleheight/2,yshifttext=-\tcboxedtitleheight/2},
	minipage boxed title=2.5cm,
	boxed title style={%
			colback=white,
			size=fbox,
			boxrule=1pt,
			boxsep=2pt,
			underlay={%
					\coordinate (dotA) at ($(interior.west) + (-0.5pt,0)$);
					\coordinate (dotB) at ($(interior.east) + (0.5pt,0)$);
					\begin{scope}
						\clip (interior.north west) rectangle ([xshift=3ex]interior.east);
						\filldraw [white, blur shadow={shadow opacity=60, shadow yshift=-.75ex}, rounded corners=2pt] (interior.north west) rectangle (interior.south east);
					\end{scope}
					\begin{scope}[gray!80!black]
						\fill (dotA) circle (2pt);
						\fill (dotB) circle (2pt);
					\end{scope}
				},
		},
	#1,
}
{prob}

\newcommand{\boldparagraph}[1]{\vspace{5pt}\noindent{}\textbf{#1.}}

\xdefinecolor{maroon}{rgb}{0.65,0.06,0.17}

\newcommand{\mb}[1]{\mathbb{#1}}
\newcommand{\mc}[1]{\mathcal{#1}}
\newcommand{\ms}[1]{\mathscr{#1}}

\newcommand{\mbf}[1]{\ensuremath{\mathbf{#1}}}

\renewcommand{\vec}[1]{\mathbf{#1}}
\newcommand{\mat}[1]{\mathbf{{#1}}}
\newcommand*{\pinv}{^{\dagger}}         
\newcommand*{\tran}{^{\mkern-1.5mu\mathsf{T}}}                
\newcommand{\sqwnorm}[2]{\left\| {#1} \right\|^2_{#2}}        

\newcommand{\R}{\mathbb{R}}

\newcommand{\e}{\varepsilon}

\DeclareMathOperator*{\argmax}{arg\,max}
\DeclareMathOperator*{\argmin}{arg\,min}

\newcommand{\state}{u}
\newcommand{\test}{p}
\newcommand{\invparam}{m}
\newcommand{\obs}{\mbf{y}}
\newcommand{\noise}{\bm{\eta}}
\newcommand{\robustparam}{\bm{\theta}}
\newcommand{\avgrobustparam}{\overline{\robustparam}}
\newcommand{\optrobustparam}{\robustparam^{\rm opt}}

\newcommand{\vzero}{\mbf{0}}

\newcommand{\design}{\bm{\xi}}
\newcommand{\scalardesign}{\xi}
\newcommand{\optdesign}{\design^{\rm opt}}
\newcommand{\alldesign}{\design^{\rm all}}
\newcommand{\sensor}{s}
\newcommand{\policy}{\mbf{p}}
\newcommand{\scalarpolicy}{p}

\newcommand{\scalarpolicyweights}{w}
\newcommand{\optpolicy}{\mbf{p}^{\rm opt}}
\newcommand{\robustoptpolicy}{\mbf{p}_{\robustparam}^{\rm opt}}

\newcommand{\baseline}{b}
\newcommand{\optbaseline}{b^{\rm opt}}

\newcommand{\hatstate}{\hat{\state}}
\newcommand{\hattest}{\hat{\test}}

\newcommand{\tildestate}{\tilde{\state}}
\newcommand{\tildetest}{\tilde{\test}}

\newcommand{\Ndata}{{\rm N_{d}}}
\newcommand{\Nrobust}{{\rm N_{\robustparam}}}
\newcommand{\Nsaa}{{\rm N_{SAA}}}
\newcommand{\Nens}{{\rm N_{ens}}}
\newcommand{\Nbudget}{{\rm N_{b}}}

\DeclareMathOperator{\diag}{diag}

\DeclareMathOperator{\range}{range}

\newcommand{\weakpde}{a}
\newcommand{\Obs}{\bm{\mc{Q}}}
\newcommand{\POM}{\bm{\mc{F}}}
\newcommand{\JPOM}{\bm{\mc{J}}}

\newcommand{\KLD}{D_{\rm KL}}
\newcommand{\EKLD}{\overline{\KLD}}

\newcommand{\Hm}{\mc{H}_{\rm m}}
\newcommand{\ppHm}{\widetilde{\mc{H}}_{\rm m}}
\newcommand{\I}{\mc{I}}

\newcommand{\utility}{\mc{U}}
\newcommand{\stochobj}{\mathfrak{U}}
\newcommand{\proj}{P}

\newcommand{\lowrankig}{\KLD^{(r)}}
\newcommand{\lowrankeig}{\overline{\lowrankig}}

\newcommand{\Expectation}{\mathbb{E}} 
\renewcommand{\Probability}{\mathbb{P}} 
\newcommand{\Variance}{\mathbb{V}} 
\newcommand{\budget}{Z} 

\newcommand{\statespace}{\ms{U}}
\newcommand{\testspace}{\ms{V}}
\newcommand{\paramspace}{\ms{M}}
\newcommand{\dataspace}{\ms{Y}}
\newcommand{\designspace}{\{0, 1\}^{\Ndata}}
\newcommand{\relaxeddesignspace}{[0, 1]^{\Ndata}}
\newcommand{\SNB}{\mathcal{S}(\Nbudget)}
\newcommand{\sensorspace}{\mc{S}}
\newcommand{\policyspace}{[0, 1]^{\Ndata}}
\newcommand{\robustspace}{\Theta}
\newcommand{\finiterobustspace}{\overline{\robustspace}}

\newcommand{\inp}[2]{\left\langle #1, #2 \right\rangle}

\newcommand{\prior}{\pi_{\rm pr}}
\newcommand{\priormeasure}{\mu_{\rm pr}}
\newcommand{\priormean}{\invparam_{\rm pr}}
\newcommand{\priorcov}{\mc{C}_{\rm pr}}
\newcommand{\priorcovinv}{\priorcov^{-1}}
\newcommand{\priorcovsqrt}{\priorcov^{1/2}}
\newcommand{\priorcovinvsqrt}{\priorcov^{-1/2}}

\newcommand{\cameronmartin}{\mc{E}}

\newcommand{\likelihood}{\pi_{\rm like}}
\newcommand{\noisecov}{\mat{\Gamma}_{\rm n}}
\newcommand{\noisecovinv}{\noisecov^{-1}}
\newcommand{\weightednoisecov}{\widehat{\mat{\Gamma}}_{\rm n}}

\newcommand{\postmeasure}[1]{\mu_{\rm post}^{#1}}
\newcommand{\lapostmeasure}[1]{\hat{\mu}_{\rm post}^{#1}}
\newcommand{\postmean}{\invparam_{\rm post}}
\newcommand{\postcov}{\mc{C}_{\rm post}}
\newcommand{\postcovinv}{\mc{C}_{\rm post}^{-1}}

\newtoggle{arxiv}
\toggletrue{arxiv}

\begin{document}

\maketitle

\begin{abstract}
	We consider robust optimal experimental design (ROED) for nonlinear Bayesian inverse
	problems governed by partial differential equations (PDEs).
	An optimal design is one that maximizes some utility quantifying the quality of the
	solution of an inverse problem.
	However, the optimal design is dependent on elements of the inverse problem such as
	the simulation model, the prior, or the measurement error model.
	ROED aims to produce an optimal design that is aware of the additional uncertainties
	encoded in the inverse problem and remains optimal even after variations in them.
	We follow a worst-case scenario approach to develop a new framework for robust optimal
	design of nonlinear Bayesian inverse problems.
	The proposed framework
	a) is scalable and designed for infinite-dimensional Bayesian nonlinear
	inverse problems constrained by PDEs;
	b) develops efficient approximations of the utility, namely, the expected information
	gain;
	c) employs eigenvalue sensitivity techniques to develop analytical forms and efficient
	evaluation methods of the gradient of the utility with respect to the uncertainties we
	wish to be robust against; and
	d) employs a probabilistic optimization paradigm that properly defines and efficiently
	solves the resulting combinatorial max-min optimization problem.
	The effectiveness of the proposed approach is illustrated for optimal sensor
	placement problem in an inverse problem governed by an elliptic~PDE.
\end{abstract}



\section{Introduction}
\label{sec:introduction}
The Bayesian approach to inverse problems is ubiquitous in the uncertainty
quantification and computational science communities.
The quality of the solution to such inverse problems is highly dependent on the
design of the data collection mechanism.
However, data acquisition is often expensive. This can put severe limits
on the amount of data that can be collected.
Thus, it is important to allocate the limited data collection resources optimally.
This can be formulated as an optimal experimental design (\textbf{OED})
problem~\cite{
	Atkinson_Donev_1992,
	Chaloner_Verdinelli_1995,
	Cox_1992,
	Ucinski_2005,
	Federov_2010}.
An OED problem seeks to identify experiments that optimize the statistical
quality of the solution to the inverse problem. While OED may be applied to a
variety of observation configurations, in this work we focus on optimal sensor
placement.  Fast and accurate OED methods for optimal design of inverse problems
governed by partial differential equations (PDEs) has been a topic of interest
over the past couple of decades; see~\cite{Alexanderian_2021} for a review of
the literature on such methods.

Inverse problems arising from complex engineering applications typically have
misspecifications and/or uncertainties in hyperparameters defining the inverse problem.
These hyperparameters, henceforth denoted by $\robustparam$,
can have a significant impact on the quality of parameter estimation.
This has sparked interest in efforts such as~\cite{Kaipio_Kolehmainen_13,
	Kolehmainen_Tarvainen_Arridge_EtAl_11,
	Mozumder_Tarvainen_Arridge_EtAl_16}
that consider Bayesian inversion under various modeling uncertainties.
See also~\cite{
	Chowdhary_Tong_Stadler_Alexanderian_2024,
	Darges_Alexanderian_Gremaud_2023,
	Sunseri_Alexanderian_Hart_Waanders_2024%
},
which consider analyzing the sensitivity of the solution of a Bayesian inverse
problem to uncertain parameters in the prior, likelihood, and governing
equations.
The uncertainty in the hyperparameters, especially the most
influential ones, needs to be accounted for in the OED problem as well.
This can be addressed by a robust OED (\textbf{ROED}) approach.
This work develops a novel ROED approach for nonlinear Bayesian inverse problems
governed by PDEs with infinite-dimensional parameters.
For consistency, hereafter, we will use the term \emph{uncertain parameter} for
any uncertain element of the inverse problem considered in an ROED framework. On
the other hand, we use the term \emph{inversion parameter} for the parameter
being estimated in an inverse problem.

\boldparagraph{Related Work}
There have been several approaches to ROED in literature.
For example, the efforts~\cite{
	Go_Isaac_2022,
	Pronzato_Walter_1988,
	Telen_Logist_VanDerlinden_VanImpe_2012,
	Telen_Vercammen_Logist_Van_Impe_2014%
}
consider finding an optimal design against a statistical average over the distribution of
the uncertain parameters. This is also related to the approach
in~\cite{Koval_Alexanderian_Stadler_2020}, which formulates the OED problem for
parameterized linear inverse problems as an optimization problem under uncertainty
problem.
Other related efforts include~\cite{
	Alexanderian_Nicholson_Petra_2024,
	Alexanderian_Petra_Stadler_Sunseri_2021,
	Bartuska_Espath_Tempone_2022%
}.
In this work, we seek to guard against the worst-case scenarios, and thus 
adopt Wald's ``max-min'' model~\cite{Wald_1945}, which seeks a design that is
optimal against a lower bound of the objective over admissible values of
the uncertain parameter.
A generic statement of the ROED problems under study is as follows:
\begin{problem}{}{roed} 
Consider the set of candidate sensor locations $\sensorspace = \{\sensor_1, \sensor_2,
	\ldots, \sensor_{\Ndata}\}$, and let $\Nbudget \ll \Ndata$ be the budget constraint on
the number of sensors. Let $\design \in \designspace$ be a binary encoding of the
observational configuration such that $\scalardesign_i$ determines whether $s_i$ is
active, and let $\robustparam \in \robustspace$ be the uncertain parameter.
The ROED problem is defined as the optimization problem
\begin{subequations}
	\begin{equation}
		\max_{\design \in \SNB}
		\min_{\robustparam \in \robustspace} \,
		\utility(\design, \robustparam) \,,
	\end{equation}
	where
	\begin{equation} \label{eq:SNB}
		\SNB = \bigg\{
		\design \in \designspace:  \sum_{i=1}^{\Ndata} \scalardesign_i = \Nbudget
		\bigg\} \,,
	\end{equation}
\end{subequations}
and the utility (objective) 
$\utility$ is chosen to quantify the quality of the design.
\end{problem}

While conservative ROED approaches has been studied in the
past~\cite{
	Biedermann_Dette_2003,
	Dette_Melas_Pepelyshev_2003,
	Pronzato_Walter_1988,
	Rojas_Welsh_Goodwin_Feuer_2007%
}, most approaches were not designed to scale to large-scale inverse problems or large
binary design spaces.
The work~\cite{Attia_Leyffer_Munson_2023} recasts the conservative (worst-case-scenario)
max-min formulation of the ROED problem  into a probabilistic optimization framework.
A desirable aspect of this development is in providing a scalable computational framework.
The approach proposed in~\cite{Attia_Leyffer_Munson_2023}, however, focuses on ROED for
linear inverse problems.
In this article, we extend the probabilistic ROED approach in \cite{Attia_Leyffer_Munson_2023} to
nonlinear Bayesian inverse problems.

\boldparagraph{Our approach and contributions}
For the utility $\utility$, we employ the expected information gain (\textbf{EIG}).
In this context, information gain is defined by the Kullback--Leibler
divergence~\cite{Kullback_Leibler_1951} of the posterior from the prior,
\begin{subequations}\label{eqn:utility-function}
	\begin{equation} \label{eq:information-gain}
		\KLD(\postmeasure{\obs} \| \priormeasure)
		:= \int \log \left(
		\dv{\postmeasure{\obs}}{\priormeasure}
		\right)
		\dd{\postmeasure{\obs}}.
	\end{equation}
	And the EIG is defined by
	\begin{equation} 
		\EKLD := \mb{E}_{\obs} \left[
			\KLD(\postmeasure{\obs} || \priormeasure)
			\right] \, .
	\end{equation}
\end{subequations}

Here, $\priormeasure$ is the prior distribution law of the inversion parameter
and $\postmeasure{\obs}$ is the posterior measure.
This utility \eqref{eqn:utility-function} admits a
computationally tractable closed-form expression in the case of Bayesian linear inverse
problems \cite{Alexanderian_Gloor_Ghattas_2016}.
However, there is no such expression in the case of nonlinear inverse problems.

In this work, we present an approximation framework for estimating the utility function
\eqref{eqn:utility-function} for infinite-dimensional Bayesian nonlinear inverse
problems governed by PDEs.
This framework enables efficient evaluation of the objective $\utility$ and its derivative
with respect to the uncertain parameter which is essential for the proposed ROED approach.
Our formulations of the utility function and its gradient are based on low-rank approximations
and adjoint-based eigenvalue sensitivity
analysis~\cite{Wu_Chen_Ghattas_2023,Chowdhary_Tong_Stadler_Alexanderian_2024}.
To enforce the budget constraint, the method in~\cite{Attia_Leyffer_Munson_2023} leveraged
a soft-constraint by adding a penalty term to the optimization objective function.
However, as pointed out in \cite{Attia_2024} this approach requires challenging tuning of the
penalty parameter.
Our proposed ROED approach leverages ideas presented in a newly developed probabilistic approach
for budget-constrained binary optimization~\cite{Attia_2024}.
We summarize the novel contributions of this article as follows: we present an ROED framework for
infinite-dimensional nonlinear Bayesian inverse problems by
\begin{enumerate}
	\item [1)] developing a scalable 
	      framework for evaluation and differentiation
	      of the EIG \eqref{eqn:utility-function} chosen as the utility function for nonlinear ROED; and
	\item [2)] developing a budget-constrained probabilistic max-min optimization framework that does not rely on penalty methods,
	      hence eliminating the need for an expensive penalty-parameter tuning stage.
\end{enumerate}
%

\boldparagraph{Article organization}
\Cref{sec:background} provides the requisite background for infinite-dimensional
Bayesian inverse problems constrained by PDEs and ROED.
\Cref{sec:roed-for-nonlinear-inverse-problems-governed-by-pdes} presents our
proposed probabilistic approach for nonlinear ROED.
Computational results demonstrating the effectiveness of the proposed methods
are in \Cref{sec:numerical-experiments}.
Concluding remarks are outlined in \Cref{sec:conclusion}.

\section{Preliminaries}
\label{sec:background}
In this section, we present the necessary mathematical background and notation for
this work.
We review infinite-dimensional Bayesian inference in
\Cref{subsec:pde-constrained-bayesian-inverse-problems} and
ROED in \Cref{subsec:oed}.

\subsection{PDE-Constrained Bayesian Inverse Problems}
\label{subsec:pde-constrained-bayesian-inverse-problems}
Consider the following PDE model:
given \emph{inversion parameter} $\invparam \in \paramspace$,
find \emph{state} $\state \in \statespace$ such that
\begin{equation} \label{eq:weak-pde}
	\weakpde(\state, \invparam, \test)
	= 0,
	\quad \forall\, \test \in \testspace
	\, ,
\end{equation}
where $\test \in \testspace$ is referred to as the \emph{test function}.
The functional $\weakpde: \statespace \times \paramspace \times \testspace \to \R$
is commonly referred to as the weak form and is linear in
$\test$. However, as we consider nonlinear models in the present work,
it is potentially nonlinear in both $\state$ and $\invparam$.

Furthermore, we assume a data model of the following form:
\begin{equation} \label{eq:data-model}
	\obs = \Obs(\state) + \noise,
\end{equation}
where $\state$ satisfies \eqref{eq:weak-pde}, and $\Obs: \statespace \to \R^d$ is an
observation operator that maps the state variable $\state$ to the data $\obs \in \dataspace$ at
$\Ndata$ observation points.
For simplicity, we assume that each sensor observes only one prognostic variable, and thus
the dimension of the observation vector is equal to the number of sensors.
This simplification, however, does not limit the formulation nor the approaches presented in this work.
Finally, we assume $\noise \sim \mc{N}(\vzero, \noisecov)$ in~\eqref{eq:data-model}.

We also define the parameter to observable map $\POM: \paramspace \to \R^d$ as
\begin{equation}\label{eq:parameter-to-observable}
	\POM(\invparam) = \Obs(\state(\invparam)) \,,
\end{equation}
which maps the inversion parameter $\invparam$ onto the observation space.
Henceforth, we assume that $\POM$ is Fr\'echet differentiable with respect to
$\invparam$.
To formulate the posterior law of the inversion parameter, we assume a Gaussian
prior measure $\priormeasure = \mc{N}(\priormean, \priorcov)$ with
$\priorcov: \paramspace \to \paramspace$ a strictly positive self-adjoint
operator of trace class; see, e.g., \cite{Stuart_2010} for details.
Here, $\priormean \in \cameronmartin$ where $\cameronmartin =
	\range(\priorcovsqrt)$ is a Cameron-Martin space, induced by the prior measure,
and is equipped with the following inner product
\begin{equation}
	\label{eq:cameron-martin-inner-product}
	\inp{x}{y}_{\priorcovinv}
	= \inp{\priorcovinvsqrt x}{\priorcovinvsqrt y}_{\!\!\paramspace}
	, \quad \forall x,y \in \cameronmartin \,.
\end{equation}

These assumptions on the data model and prior, along with Bayes' rule, define a
posterior measure $\postmeasure{\obs}$ on $\paramspace$ given by the
Radon--Nikodym derivative
\begin{equation}\label{eq:bayes}
	\dv{\postmeasure{\obs}}{\priormeasure}
	\propto \likelihood(\obs | \invparam).
\end{equation}
Note that under the additive Gaussian noise model \eqref{eq:data-model}, the likelihood
$\likelihood$ satisfies
\begin{equation} \label{eq:likelihood}
	\likelihood(\obs | \invparam)
	\propto
	\exp \left(
	-\frac{1}{2} \norm{\obs - \POM(\invparam)}_{\noisecovinv}^2
	\right),
\end{equation}
where we have used the weighted norm $\sqwnorm{\vec{x}}{\noisecovinv} :=
	\vec{x}\tran\noisecovinv\vec{x}$.
If the parameter-to-observable map $\POM$ is linear, then it can be demonstrated
\cite{Stuart_2010} that the posterior measure is Gaussian $\postmeasure{\obs} =
	\mc{N}(\postmean, \postcov)$ with
\begin{equation} \label{eq:linear-posterior}
	\postcov
	=
	\left( \priorcovinv + \POM^* \noisecovinv \POM \right)^{-1}  \,;
	\qquad 
	\postmean
	=
	\postcov
	\big(
	\POM^* \noisecovinv \obs + \priorcovinv \priormean
	\big) \,,
\end{equation}
where $\POM^{*}$ is the adjoint of $\POM$.
In this work, however, we explicitly consider the case where the
parameter-to-observable map $\POM$ is nonlinear.
In this case, $\postmeasure{\obs}$ is generally not available in closed form,
and must be approximated through various methods such as Markov Chain Monte
Carlo (MCMC) or Variational Inference (VI).

Recall, our objective, as described in \Cref{prob:roed}, is to determine an
optimal sensor configuration that's robust to misspecified hyperparameters.
At this point, we have not been explicit regarding the parameterization
of said hyperparameters $\robustparam$.
The formulation of ROED permits $\robustparam$ to be anywhere in the elements
of the inverse problem.
For example, $\robustparam$ may parameterize the prior distribution $\priormeasure$,
the likelihood $\likelihood$, or even parameter-to-observable map $\POM$ itself.
Only a single contribution of this paper requires knowledge of the exact
parameterization, see
\Cref{subsubsec:differentiating-utility-via-variational-tools}.
Therefore, until then, we assume $\postmeasure{\obs} = \postmeasure{\obs}(\robustparam)$.

\subsection{Robust Optimal Experimental Design}
\label{subsec:oed}
As stated in \Cref{prob:roed}, our goal is to find a robust optimal design
$\optdesign$ that solves the ROED optimization problem
\begin{equation}
	\label{eq:roed}
	\max_{\design \in \SNB}
	\min_{\robustparam \in \robustspace} \,
	\utility(\design, \robustparam) \,,
\end{equation}
where $\utility$ is the chosen utility function.
Let us contrast this with the traditional (non-robust) optimal design that
solves the binary optimization problem
\begin{equation}
	\label{eq:oed}
	\max_{\design \in \SNB} \utility(\design).
\end{equation}
The ROED problem \eqref{eq:roed} can be viewed as a bi-level optimization problem, with
the outer optimization layer bearing resemblance to the traditional OED problem
\eqref{eq:oed}.
Hence, solving the ROED problem is considerably more challenging than
the traditional OED problem.

Some of the commonly used techniques for solving
traditional OED problems are not suitable for the ROED problem.
Notably, we recall that a commonly used technique for solving \eqref{eq:oed} is a
relaxation approach. In this approach
the design space is relaxed from a binary space $\designspace$ to a
continuum $\relaxeddesignspace$.
This enables the use of gradient-based optimization techniques.  However, as
noted in~\cite{Attia_Leyffer_Munson_2023}, a naive reformulation of
\eqref{eq:roed} into this relaxed form is incorrect because
\[
	\argmax_{\design \in \designspace}
	\min_{\robustparam \in \robustspace} \,
	\utility(\design, \robustparam)
	\not\subset
	\argmax_{\design \in \relaxeddesignspace}
	\min_{\robustparam \in \robustspace} \,
	\utility(\design, \robustparam).
\]
Specifically, it was demonstrated in \cite{Attia_Leyffer_Munson_2023} that
relaxation of the binary design in \eqref{eq:roed} results in a different optimization problem
with different optimal set than the solution of \eqref{eq:roed}.
To overcome this challenge, \cite{Attia_Leyffer_Munson_2023} introduced an extension
of a stochastic optimization framework \cite{Attia_Leyffer_Munson_2022} to the ROED setting and
demonstrated that the resulting probabilistic ROED formulation is equivalent to
the original binary max-min ROED optimization problem \eqref{eq:roed}.
\Cref{prob:old-stochastic-roed} summarizes this probabilistic ROED formulation.
In the present work, we consider extensions of such formulations
to budget-constrained ROED for nonlinear Bayesian inverse problems.
\begin{problem}{}{old-stochastic-roed} 
Consider the set of candidate sensor locations
$\sensorspace = \{\sensor_1, \sensor_2, \ldots, \sensor_{\Ndata}\}$. Let
$\design \in \designspace$ be a binary encoding of the
observational configuration such that $\scalardesign_i$ determines if $s_i$ is
active. Let $\robustparam \in \robustspace$ be the uncertain parameter we seek to be
robust against.
The probabilistic ROED approach views $\design$ as a random variable endowed with
a multivariate Bernoulli distribution
$\Probability(\design | \policy)$
parameterized by $\policy \in \policyspace$\,,
where $\scalarpolicy_i \in [0,1]$ is the probability of activating
$s_i$.

The probabilistic ROED problem aims to find a policy $\optpolicy$
that solves
\begin{equation} \label{eq:robust-oed}
	\max_{\policy \in \policyspace}
	\stochobj(\policy)
	:=
	\Expectation_{\design \sim \Probability(\design | \policy)} \left[
		\min_{\robustparam \in \robustspace} \,
		\utility(\design, \robustparam)
		\right] \,.
\end{equation}
Here, $\Probability(\design | \optpolicy)$ yields the solution of the binary ROED
problem~\eqref{eq:roed}.
\end{problem}

The algorithmic approach presented in \cite{Attia_Leyffer_Munson_2023} for solving
\eqref{eq:robust-oed} relies on an efficient sampling based approach that was originally
introduced in \cite{Levitin_Polyak_1966}.
We defer the majority of the details of this approach to
\cite{Attia_Leyffer_Munson_2023,Levitin_Polyak_1966}.
For clarity, and to highlight the contributions of this work,
we only provide a brief overview of the algorithm in this section.
A complete algorithmic statement of our proposed approach which extends this
sampling-based approach is provided in \Cref{subsec:computational-complexity}.

The sampling based approach for solving \eqref{eq:robust-oed} is an iterative procedure
that alternates between solving an outer optimization over the policy $\policy$,
and an inner optimization problem over the uncertain parameter $\robustparam$.
In the outer optimization stage, the expectation is approximated by using a finite set
of samples from $\finiterobustspace\subset\robustspace$.
The sample $\finiterobustspace$ is then expanded by solving the inner optimization problem
of the uncertain parameter.
Thus, at iteration $k$ of the optimization procedure, with the finite sample of the
uncertain parameter $\finiterobustspace^{(k)}\subset\robustspace$,
the outer optimization problem seeks a $\policy^{(k)}$ that
maximizes
\begin{equation}
	\stochobj^{(k)}(\policy)
	=
	\Expectation_{\design \sim \Probability(\design | \policy)}
	\left[
		\min_{\robustparam \in \finiterobustspace^{(k)}} \,
		\utility(\design, \robustparam)
		\right] \,.
\end{equation}
At the same iteration $k$, the inner optimization problem seeks a $\robustparam^{(k)}$
that minimizes $\utility$ over $\robustspace$ using designs sampled from $\Probability(\design |
	\policy^{(k)})$; this minimizer is added to the set $\finiterobustspace^{(k)}$.
The algorithm follows a gradient-based approach for solving both the outer and the inner optimization problems.
For the outer optimization problem, a stochastic gradient is used which requires the gradient of the
probability model with respect to its parameter $\nabla_{\policy} \Probability(\design | \policy)$.
The inner optimization problem, however, requires the gradient of $\utility$ with respect
to the uncertain parameter $\robustparam$.
In this approach, however, the utility function $\utility$ must be differentiable
with respect to the uncertain parameter.

With that in mind, we highlight a few critical benefits and limitations of the
probabilistic ROED approach defined by \Cref{prob:old-stochastic-roed}.
A major \emph{advantage} of \Cref{prob:old-stochastic-roed} demonstrated in
\cite{Attia_Leyffer_Munson_2023} is its scalability with respect to $\Ndata$ and
$\Nrobust$. Additionally, this approach opens the way for use of gradient-based
optimization methods in the outer design optimization stage without requiring
derivatives of the utility function with respect to $\design$.
A key \emph{limitation} of
this approach in its present formulation, is that the distribution
$\Probability(\design | \policy)$ does not impose any budget constraint on the number of active sensors.
Hence, any budget constraint would typically be enforced through a penalty term in the utility function.
This necessitates an expensive hyperparameter tuning phase.
Likewise, during the optimization procedure, designs sampled from the distribution
$\Probability(\design | \policy)$ are not guaranteed to satisfy the budget constraint, hence
potentially spending computational resources on infeasible designs.
Finally, as mentioned earlier, although gradients of $\utility$ with respect
to $\design$ are not required, gradients of the utility function $\utility$ with respect
to the uncertain parameter $\robustparam$ are required for the inner optimization stage.
Overcoming these challenges for ROED for Bayesian nonlinear inverse problems is the
primary objective of the contributions of this work.

\section{Robust Optimal Experimental Design for Bayesian Nonlinear Inverse Problems}
\label{sec:roed-for-nonlinear-inverse-problems-governed-by-pdes}
In this section we propose a scalable ROED approach
for nonlinear infinite-dimensional inverse problems under budget-constraints.
This begins with a new formulation of probabilistic ROED optimization problem
with a budget-constrained probability distribution in
\Cref{subsec:budget-constrained-stochastic-roed}.
This formulation of the ROED optimization problem is applicable to any choice
of the utility function.
Then, in \Cref{subsec:approximating-the-expected-information-gain},
we focus on EIG as the utility function.  In that section, we
discuss an approximation framework to enable fast evaluation
(\Cref{subsubsec:low-rank-and-fixed-map-approximation}) and differentiation
(\Cref{subsubsec:differentiating-utility-via-variational-tools}) with respect to
the uncertain parameters of the EIG.

\subsection{Budget-Constrained Stochastic Robust OED}
\label{subsec:budget-constrained-stochastic-roed}
In this section, we introduce a new formulation for ROED that enforces a budget
constraint on the number of active sensors.
To do so, we first introduce the conditional Bernoulli model developed in \cite{Attia_2024}.
We restate its definition in our notation:
\begin{definition}
	\label{def:conditional-bernoulli}
	Let $\design = (\scalardesign_1, \dots, \scalardesign_{\Ndata}) \in
		\designspace$ be a multivariate Bernoulli random variable parameterized by
	the policy $\policy = (\scalarpolicy_1, \dots, \scalarpolicy_{\Ndata}) \in
		\policyspace$.
	Let $\budget \equiv \budget(\design) = \sum_{i=1}^{\Ndata}
		\scalardesign_i$ be the total number of active (equal to $1$) entries in $\design$,
	and define
	\begin{equation}
		S = \{1, \dots, \Ndata\};
		\ 
		O = \{i \in S: \scalarpolicy_i \!=\! 0\};
		\ 
		I = \{i \in S: \scalarpolicy_i \!=\! 1\};
		\, 
		T = S \setminus \{O \cup I\} \,.
	\end{equation}
	Then, the probability mass function (PMF) of the conditional
	Bernoulli model is:
	\begin{equation}
		\label{eq:bernoulli-budget-pmf}
		\Probability(\design | \policy, \budget = z)
		=
		\begin{cases}
			\frac{
				\prod\limits_{i \in T} \scalarpolicyweights^{\scalardesign_i}
			}
			{ R(z - |I|, T) }
			,  &
			\textrm{if }
			\scalardesign_j = \scalarpolicy_j, \forall j \in \{I \cup O\}
			\ \text{and} \
			\sum\limits_{j \in T} \scalardesign_j = z - |I|
			\\
			0, & \textrm{otherwise}
		\end{cases}
	\end{equation}
	where
	\begin{equation}
		\label{eq:r-function}
		R(k, A)
		=
		\sum_{\substack{B \subseteq A \\ |B| = k}}
		\prod_{i \in B} \scalarpolicyweights_i
		\,;  \qquad 
		\scalarpolicyweights_i
		=
		\frac{\scalarpolicy_i}{1 - \scalarpolicy_i}
		, \forall i \in \{1, \dots, \Ndata\} \,.
	\end{equation}
\end{definition}

Details regarding the fast evaluation and differentiation of the PMF in
\Cref{def:conditional-bernoulli} can be found in~\cite{Attia_2024}.  Using this
conditional distribution, we construct a modification of the probabilistic
robust OED framework discussed in \Cref{prob:old-stochastic-roed}.  Our
proposed problem formulation, stated in \Cref{prob:budget-constrained-stochastic-roed},
enforces the budget constraint \eqref{eq:SNB} without the need for a penalty
term in the utility function.

\begin{problem}{}{budget-constrained-stochastic-roed} 
Let $\sensorspace = \{\sensor_1, \ldots, \sensor_{\Ndata}\}$ be the set of
candidate sensor locations, $\Nbudget \ll \Ndata$ be the budget constraint, $\design
	\in \designspace$ be a binary encoding of the observational configuration such that
$\scalardesign_i$ determines if $s_i$ is active, and $\robustparam \in
	\robustspace$ be the uncertain parameter we seek to be robust against.

Now, let us assume that $\design$ is a random variable endowed with the conditional
Bernoulli distribution $\Probability(\design | \policy, \budget = \Nbudget)$
as defined in \Cref{def:conditional-bernoulli}.
Then, the budget-constrained probabilistic ROED problem replaces \Cref{prob:roed}
with the following policy optimization problem:
\begin{equation}
	\label{eq:budget-constrained-stochastic-roed}
	\max_{\policy \in \policyspace}
	\stochobj(\policy)
	:=
	\Expectation_{\design \sim \Probability(\design | \policy, \budget = \Nbudget)}
	\left[
		\min_{\robustparam \in \robustspace} \,
		\utility(\design, \robustparam)
		\right] \,.
\end{equation}
We denote $\stochobj$ as the stochastic objective.
\end{problem}

To solve \Cref{prob:budget-constrained-stochastic-roed}, we leverage the same
sampling-based approach described in \Cref{subsec:oed}.
However, given the modifications to the probability distribution,
we need to re-derive the necessary components of
the algorithm.
That is, we need to re-derive the gradients of the stochastic objective $\stochobj$
with respect to the policy parameter $\policy$ for the outer optimization stage and the gradients
of the utility function $\utility$ with respect to the uncertain parameter $\robustparam$ for the inner optimization stage.
We defer the discussion of the latter to the next section,
as it depends on the specific form of $\utility$.

Now, let us consider the computation of the gradients of $\stochobj$ with respect to
$\policy$. Note, by the definition of $\stochobj$
\begin{equation}\label{eq:exact-roed-gradient}
	\grad_{\policy} \stochobj(\policy)
	=
	\Expectation_{\design \sim \Probability(\design | \policy, \budget = \Nbudget)}
	\left[
		\min_{\robustparam \in \robustspace} \,
		\utility(\design, \robustparam) \,\,
		\grad_{\policy} \log \Probability(\design | \policy, \budget = \Nbudget)
		\right],
\end{equation}
where we have leveraged the equality $\grad_\policy f(\policy) = f(\policy)
	\grad_\policy \log f(\policy)$, for $\policy \in (0, 1)^{\Ndata}$.
Note that if $\scalarpolicy_i \in\{0, 1\}$, the corresponding partial derivative
is set to zero; see~\cite{Attia_Leyffer_Munson_2022}.

Now, to evaluate \cref{eq:exact-roed-gradient} directly would be intractable as the
underlying discrete space is of cardinality $\binom{\Ndata}{\Nbudget}$. Instead, we
adopt a stochastic gradient approximation approach. Specifically, given samples
\(
\{
\design[k] \sim \Probability( \design | \policy, \budget = \Nbudget )
| k = 1, \dots, \Nens
\}
\), the stochastic approximation of the gradient \eqref{eq:exact-roed-gradient} is given by
\begin{equation}
	\label{eq:stochastic-roed-gradient}
	\grad_{\policy} \stochobj(\policy)
	\approx
	\frac{1}{\Nens}
	\sum_{k=1}^{\Nens}
	\left[
		\min_{\robustparam \in \robustspace} \,
		\utility(\design[k], \robustparam) \,\,
		\grad_{\policy} \log \Probability(\design[k] | \policy, \budget = \Nbudget)
		\right].
\end{equation}

Finally, as seen in \cite{Attia_Leyffer_Munson_2022}, the performance of stochastic
gradient estimator is greatly enhanced by using variance reduction techniques
such as optimal baseline.
This technique replaces the utility function
$\utility(\design, \robustparam)$ with $\utility(\design, \robustparam) - \baseline$,
where baseline $\baseline$ is a constant scalar selected to minimize the variance of the gradient
estimator.
Towards determining an optimal value of that baseline,
we define the stochastic objective with a
baseline and its policy gradient as
\begin{subequations}
	\begin{align}
		\label{eq:stochastic-roed-objective-baseline}
		\stochobj^{\baseline}(\policy)
		 & =
		\Expectation_{\design \sim \Probability(\design | \policy, \budget = \Nbudget)}
		\left[
			\min_{\robustparam \in \robustspace} \,
			\left[ \utility(\design, \robustparam) - \baseline \right]
		\right],   \\
		\label{eq:stochastic-roed-gradient-baseline}
		\grad_{\policy} \stochobj^{\baseline}(\policy)
		 & \approx
		\frac{1}{\Nens}
		\sum_{k=1}^{\Nens}
		\left[
			\min_{\robustparam \in \robustspace} \,
			\left[ \utility(\design[k], \robustparam) - \baseline \right]
			\grad_{\policy} \log \Probability(\design[k] | \policy, \budget = \Nbudget)
			\right].
	\end{align}

	The optimal baseline $\optbaseline$ is found by minimizing the variance of the gradient
	estimator with respect to $\baseline$, as demonstrated in \cite{Attia_Leyffer_Munson_2023,Attia_2024}.
	Here, $\optbaseline = \max(0, \baseline^*)$ where
	\begin{equation}
		\label{eq:stochastic-budget-roed-optimal-baseline}
		\baseline^*
		=
		\frac{
			\sum\limits_{i=1}^{\Nens}
			\sum\limits_{j=1}^{\Nens}
			\utility(\design[i], \robustparam^*[i])
			\inp{
				\grad_{\policy} \log \Probability(\design[i] | \policy, \budget = \Nbudget)
			}
			{
				\grad_{\policy} \log \Probability(\design[j] | \policy, \budget = \Nbudget)
			}
		}
		{
			\Nens
			\sum\limits_{i=1}^{\Ndata}
			\frac{(1+\scalarpolicyweights_i)^4}{\scalarpolicyweights_i^2}
			(\pi_i - \pi_i^2)
		}
	\end{equation}
	and
	\begin{equation}
		\label{eq:stochastic-budget-roed-optimal-baseline-robustparam}
		\robustparam^*[i]
		=
		\argmin_{\robustparam \in \robustspace}
		\utility(\design[i], \robustparam)
		\, ; \qquad
		\pi_i
		=
		\scalarpolicyweights_1
		\frac{R(\Nbudget - 1, S \setminus \{i\})}{R(\Nbudget, S )} \,.
	\end{equation}
\end{subequations}

Naturally, the minimization problem
\eqref{eq:stochastic-budget-roed-optimal-baseline-robustparam} would be replaced by one
over $\finiterobustspace^{(k)}$ when computing the optimal baseline in the context of
the outer optimization stage of the ROED algorithm.
Details regarding the evaluation of the total variance of the gradient
$ \Variance\left[ \grad_{\policy} \log \Probability(\design | \policy, \budget = \Nbudget)\right] $
in \eqref{eq:stochastic-roed-gradient-baseline}
may be found in
\cite[Section 3.2]{Attia_2024}.
A complete algorithmic description of the solution process is given
by \Cref{alg:stochastic-budget-roed}.

\subsection{The Utility Function: Expected Information Gain}
\label{subsec:approximating-the-expected-information-gain}
Now we turn our attention to the fast estimation and differentiation of the utility
function $\utility$.
Up until now, our methods have been agnostic to the choice of utility and we
emphasize that the framework presented in
\Cref{subsec:budget-constrained-stochastic-roed} would be applicable to
any utility function.
However, henceforth we focus on an approximation framework to accelerate the
Expected Information Gain (EIG) for the purposes of ROED.

\subsubsection{Expected Information Gain}
\label{subsubsec:expected-information-gain}
The EIG is a widely used information-theoretic utility function for nonlinear
optimal experimental design \cite{Chaloner_Verdinelli_1995}.
In the present setting, the EIG is given by
\begin{equation} \label{eq:expected-information-gain}
	\EKLD
	=
	\mb{E}_{\obs}\! \left[
		\KLD(\postmeasure{\obs} \,\|\, \priormeasure)
		\right]
	= \!
	\int\limits_\paramspace
	\int\limits_\dataspace
	\KLD(\postmeasure{\obs} \,\|\, \priormeasure)
	\likelihood(\obs | \invparam)
	\dd{\obs}
	\dd{\priormeasure(m)} \,.
\end{equation}
In the case of a linear parameter-to-observable map $\POM$, $\EKLD$ attains the
following closed form expression (see; e.g.,
\cite{Alexanderian_Gloor_Ghattas_2016}):
\begin{equation}
	\label{eq:linear-expected-information-gain}
	\EKLD = \frac{1}{2} \log\det\left( \I + \ppHm \right) \,,
\end{equation}
where $\I$ is the identity operator, and $\ppHm = \priorcovsqrt \Hm
\priorcovsqrt$ is the prior preconditioned data-misfit Hessian.
This fact has been employed for both the fast evaluation
\cite{Alexanderian_Saibaba_2018} of the EIG and scalable differentiation of it
with respect to model hyperparameters
\cite{Chowdhary_Tong_Stadler_Alexanderian_2024}.
The original probabilistic approach in \cite{Attia_Leyffer_Munson_2023} has also
formulated derivatives of this expression with respect to the uncertain
parameter $\robustparam$.
In the nonlinear setting considered in this work, however, no such closed form
expression for the EIG exists, and we must proceed from the double integral in
\eqref{eq:expected-information-gain}.

Evaluating \eqref{eq:expected-information-gain} following a Monte-Carlo estimation
approach has been studied for lower dimensional problems; see, e.g.,
\cite{Rainforth_Foster_Ivanova_BickfordSmith_2024}.
%
This approach, however, does not scale well to high dimensional problems and is
thus not suitable for infinite-dimensional settings.
To proceed, we leverage the fact that the information gain from the prior to
the \emph{Laplace approximation} to the posterior, $\lapostmeasure{\obs}$, attains a
closed form expression~\cite{Chowdhary_2025}.

\boldparagraph{Laplace Approximation}
The Laplace approximation aims to approximate the posterior by an appropriate
Gaussian distribution, and is a commonly used tool in large-scale nonlinear
Bayesian inverse problems
\cite{Stuart_2010,Buithanh_Ghattas_Martin_Stadler_2013},
It is the Gaussian $\lapostmeasure{\obs} = \mc{N}(\postmean, \postcov)$, where
$\postmean$ is given by the maximum a posteriori (MAP) point
\begin{subequations}\label{eq:laplace-approximation}
	\begin{equation} \label{eq:postmean}
		\argmin_{\invparam \in \cameronmartin}
		\Phi(\invparam)
		\coloneqq
		\frac{1}{2} \norm{\obs - \POM(\invparam)}_{\noisecovinv}^2
		+ \frac{1}{2} \norm{\invparam - \priormean}_{\priorcovinv}^2 \,,
	\end{equation}
	and covariance $\postcov$ given by the inverse Hessian of
	$\Phi$ evaluated at $\postmean$
	\begin{equation} \label{eq:postcov}
		\postcovinv = \Hm(\postmean) + \priorcovinv.
	\end{equation}
	Here, $\Hm(\postmean)$ is the Hessian of the data misfit term $\frac{1}{2}
		\| \obs - \POM(\invparam) \|_{\noisecovinv}^2$ evaluated at $\postmean$ and
	has been separated out for notational convenience.
\end{subequations}

In practice, a commonly used approximation to this data-misfit Hessian is the
Gauss--Newton approximation given by
\begin{equation} \label{eq:Hm}
	\Hm(\postmean) = \JPOM^*(\postmean) \noisecovinv \JPOM(\postmean),
\end{equation}
where $\JPOM(\postmean)$ is the Jacobian of $\POM(\invparam)$ with respect to
$\invparam$ evaluated at $\postmean$.
We employ this approximation in our work, and henceforth we use the simplified notation
$\Hm$ to refer to the expression in \eqref{eq:Hm}.
We note that with the Gauss--Newton Hessian approximation \eqref{eq:Hm}, the
Laplace approximation is equivalent to the posterior measure obtained after a
linearization of the parameter-to-observable map $\POM$ at the MAP point
$\postmean$; see e.g., \cite{Alexanderian_Nicholson_Petra_2024, Chowdhary_2025}.
Note, a computationally efficient way of evaluating the action of $\Hm$ can be
derived through an adjoint-based gradient approach \cite{Plessix_2006, Bangerth_2008}.
In particular,
\begin{subequations} \label{eq:gn-data-misfit-hessian-system}
	\begin{equation} \label{eq:gn-data-misfit-hessian}
		\mc{H}(\invparam)(\hat{\invparam}, \tilde{\invparam})
		= \inp{\tilde{\invparam}}{
			\weakpde_{\invparam\invparam}(\state, \invparam, \test)
			+ \weakpde_{\invparam\state}(\state, \invparam, \test) \hat{\state}
			+ \weakpde_{\invparam\test}(\state, \invparam, \test) \hat{\test}
		},
	\end{equation}
	where
	for all $\tilde{\test} \in \testspace$
	and $\tilde{\state} \in \testspace$,
	\begin{align}
		\label{eq:gn-incr-state}
		\inp{\tilde{\test}}{\weakpde_{\test\state}(\state, \invparam, \test) \hat{\state}}
		+ \inp{\tilde{\test}}{\weakpde_{\test\invparam}(\state, \invparam, \test)
			\hat{\invparam}}
		 & = 0 \,, 
		\\
		\label{eq:gn-incr-adjoint}
		\inp{\tilde{\state}}{
			\weakpde_{\state\test}(\state, \invparam, \test) \hat{\test}
		}
		+ \inp{\tilde{\state}}{\Obs^* \noisecovinv \Obs \hat{\state}}
		 & = 0 \,. 
	\end{align}
\end{subequations}
We assume the inner products above are defined appropriately
for the underlying PDE and $\testspace$.
We provide a detailed derivation for \Cref{eq:gn-data-misfit-hessian-system} in
\Cref{sec:variational-tools}.

With this in hand, we now note that \cite[Proposition 3.3.1]{Chowdhary_2025}
\begin{subequations}\label{eq:expected-information-gain-closed-average}
	\begin{multline} \label{eq:closed-form-information-gain}
		\KLD(\lapostmeasure{\obs} \,\|\, \priormeasure)
		= \\ \frac{1}{2} \Big[
			\log\det\Big( \I + \ppHm \Big)
			- \trace \Big( \ppHm \big[ \I + \ppHm \big]^{-1} \Big)
			+ \norm{\postmean - \priormean}_{\priorcovinv}^2
			\Big] \,.
	\end{multline}
	This enables approximating the EIG \eqref{eq:expected-information-gain} by
	the sample average approximation
	\begin{equation} \label{eq:saa-eig}
		\EKLD
		\approx
		\frac{1}{\Nsaa} \sum_{i=1}^{\Nsaa}
		\KLD(\lapostmeasure{\obs_i} \,\|\, \priormeasure) \,,
	\end{equation}
	where for every $i \in \{1, \ldots, \Nsaa\}$, the data $\obs_i$ are drawn from the model
	\begin{equation}
		\obs_i = \POM(\invparam_i) + \noise_i
		\, ,
	\end{equation}
	where $\invparam_i \sim \priormeasure$ and $\noise_i \sim \mc{N}(\vzero, \noisecov)$.
\end{subequations}
This approximation has proven effective in the OED context in literature
\cite{Wu_Chen_Ghattas_2023}.

However, there remain challenges in scaling this approach to infinite-dimensional ROED.
While data-parallel, this approximation requires $\Nsaa$ MAP estimations.
Additionally, computing the first two terms of \eqref{eq:closed-form-information-gain}
involve estimating the log-determinant and trace of high-dimensional operators.
For even problems at a moderate scale, this approach may prove computationally
challenging.
Likewise, a scalable procedure for differentiating this expression with respect to
$\robustparam$ is not immediately clear.
Our final contribution is the development of additional techniques to further
approximate the EIG.
Furthermore, we introduce an adjoint-based eigenvalue sensitivity approach to
differentiating it with respect to $\robustparam$.
These two novel components, along with the discussion in
\Cref{subsec:budget-constrained-stochastic-roed}, will then be used to develop a
complete algorithmic statement of our proposed ROED approach.

\subsubsection{Low-Rank and Fixed MAP Approximation}
\label{subsubsec:low-rank-and-fixed-map-approximation}
The prior preconditioned data-misfit Hessian $\ppHm\equiv \ppHm(\design, \robustparam)$
is often low-rank.
We leverage this structure to find computationally efficient approximations of the first
two terms of the information gain \eqref{eq:closed-form-information-gain} in terms of
the dominant eigenvalues of $\ppHm$.
This type of approximation for the information gain has been utilized in prior works
such as~\cite{Alexanderian_Saibaba_2018,Chowdhary_Tong_Stadler_Alexanderian_2024,
	Wu_Chen_Ghattas_2023}.

Due to the structure of the Gauss-Newton Hessian \eqref{eq:Hm}, the prior preconditioned
data-misfit Hessian has rank of at most $\Ndata$, corresponding to the design with all
sensors active.
However, during the optimization process defined in
\Cref{prob:budget-constrained-stochastic-roed}, only designs with $\Nbudget$ active
sensors are considered. Hence, the rank of $\ppHm$ is at most $\Nbudget$.
Therefore, leveraging a randomized method with $\mc{O}(\Nbudget)$ applications of the
Hessian we can obtain the low-rank approximation
\begin{equation}
	\ppHm \phi
	= \sum_{n=1}^\infty \lambda_n \inp{\phi}{\omega_n} \omega_n
	\approx \sum_{n=1}^r \lambda_n \inp{\phi}{\omega_n} \omega_n
	, \quad \phi \in \paramspace \,,
\end{equation}
where $r \leq \Nbudget$ is some appropriately chosen integer such that $(\lambda_n,
	\omega_n)$ are the $r$ dominant eigenpairs of $\ppHm$.
This is given by the eigenproblem
\begin{equation}
	\inp{\phi}{\ppHm \omega_n} = \lambda_n \inp{\phi}{\omega_n},
	\quad \text{with} \quad
	\inp{\omega_n}{\omega_n} = 1,
	\quad \forall \phi \in \paramspace, i \in \{1, \ldots, r\} \,.
\end{equation}

Note that the eigenvalues of $\ppHm$ are dependent on the data realization $\obs$ used in the
inverse problem as well as the design $\design$ and uncertain parameter $\robustparam$.
To be precise, we make these dependencies explicit.
Namely, for data $\obs_i$ we denote the resulting MAP point by
$\postmean^i(\design, \robustparam)$ and the
prior-preconditioned data-misfit Hessian by $\ppHm^i(\design, \robustparam)$.
Likewise, we denote the dominant eigenvalues of $\ppHm^i(\design, \robustparam)$ as
$\{\lambda_n^i(\design, \robustparam)\}_{n=1}^r$.
Thus, from \eqref{eq:expected-information-gain-closed-average} it follows that
a Laplace approximation approach yields the following information gain estimate
\begin{equation*}
	\begin{aligned}
		\KLD(\obs_i, \design, \robustparam)
		\approx
		 &
		\frac{1}{2}
		\left[
			\log\det \Big( \I + \ppHm^i(\design, \robustparam) \Big)
			- \trace \Big(
			\ppHm^i(\design, \robustparam)
			\big[ \I + \ppHm^i(\design, \robustparam) \big]^{-1}
			\Big)
			\right]
		\\
		+
		 &
		\frac{1}{2}\norm{\postmean^i(\design, \robustparam) - \priormean}_{\priorcovinv}^2
		\, .
	\end{aligned}
\end{equation*}

Using the dominant eigenvalues of $\ppHm^i(\design, \robustparam)$, we can approximate
the first two terms to define the \textit{low-rank information gain} as
\begin{subequations}\label{eq:Laplace-low-rank-EIG}
	\begin{equation}
		\label{eq:low-rank-information-gain}
		\begin{aligned}
			\lowrankig(\obs_i, \design, \robustparam)
			=
			 &
			\frac{1}{2}
			\sum_{n=1}^r
			\Bigg[
				\log(1 + \lambda_n^i(\design, \robustparam))
				- \frac{\lambda_n^i(\design, \robustparam)}{1+\lambda_n^i(\design, \robustparam)}
				\Bigg]
			\\
			+
			 &
			\frac{1}{2}
			\norm{\postmean^i(\design, \robustparam) - \priormean}_{\priorcovinv}^2
			\,.
		\end{aligned}
	\end{equation}
	Thus, the \textit{low-rank EIG} is given by
	\begin{equation}
		\label{eq:saa-low-rank-eig}
		\lowrankeig(\design, \robustparam)
		=
		\frac{1}{\Nsaa} \sum_{i=1}^{\Nsaa} \lowrankig(\obs_i, \design, \robustparam) \,.
	\end{equation}
\end{subequations}

While the approximation \eqref{eq:low-rank-information-gain}
provides an efficient method for estimating the first two terms of the
information gain, it still requires a MAP point estimation $\postmean^i(\design, \robustparam)$.
In the sample average approximation for the EIG \eqref{eq:saa-low-rank-eig}, one would
therefore need to compute $\Nsaa$ MAP point solvers per evaluation, which is
computationally challenging.
In \cite{Wu_Chen_Ghattas_2023}, the authors proposed a fixed MAP point
approximation to alleviate this burden.
Let $\alldesign$ be the design with all sensors active.
Then, the fixed MAP point approximation replaces $\postmean^i(\design, \robustparam)$
with $\postmean^i(\alldesign, \robustparam)$ for every $i \in \{1, \ldots, \Nsaa\}$.
In ROED settings, however, the MAP point is also dependent on
$\robustparam$. Hence, simply fixing a nominal value for the design
parameter does not resolve the need to perform $\Nsaa$ MAP estimations per evaluation.
Therefore, we extend the fixed MAP point approximation to $\robustparam$ as well.
Typically, only have access to the uncertain parameter space $\robustspace$ through a
finite sample $\finiterobustspace$. Hence, in the present work, we propose to
additionally fix the uncertain parameter at the ensemble average of the finite sample,
hence, define the fixed MAP estimate as
\begin{equation}
	\label{eq:average-robust-parameter}
	\postmean^i(\alldesign, \avgrobustparam)\,; \qquad
	\avgrobustparam = \frac{1}{|\finiterobustspace|}\sum_{\robustparam \in \finiterobustspace} \robustparam \,.
\end{equation}
This leads to the following ROED utility function defined
using the \textit{low-rank EIG with a fixed MAP approximation}:
\begin{subequations}\label{eqn:roed_utility_function}
	\begin{align}
		\label{eq:saa-low-rank-fixed-map-eig}
		\utility(\design, \robustparam)
		 & =
		\frac{1}{\Nsaa} \sum_{i=1}^{\Nsaa}
		\hat{\utility}(\obs_i, \design, \robustparam) \,,
		\\
		\label{eq:low-rank-fixed-map-ig}
		\hat{\utility}(\obs_i, \design, \robustparam)
		 & =
		\frac{1}{2}\sum_{n=1}^r
		\left[
			\log(1 + \lambda_n^i(\design, \robustparam))
			- \frac{\lambda_n^i(\design, \robustparam)}{1+\lambda_n^i(\design, \robustparam)}
			\right]
		+ C_i\,,
	\end{align}
\end{subequations}
where $C_i =
	\frac{1}{2}\|\postmean^i(\alldesign, \avgrobustparam) - \priormean\|_{\priorcovinv}^2$.
We emphasize that the fixed MAP points $\{ \postmean^i(\alldesign, \avgrobustparam) : i
	= 1, \ldots, \Nsaa\}$ are computed offline. That is to say, they're computed once at
the beginning of the computation and are reused for subsequent evaluations of the
utility function.
Thus, to evaluate the utility function \eqref{eqn:roed_utility_function} across
different values of $\design$ and $\robustparam$, only the randomized eigendecomposition
need be performed $\Nsaa$ times.

\subsubsection{Differentiation via Variational Tools}
\label{subsubsec:differentiating-utility-via-variational-tools}
Finally, for the inner optimization of the stochastic ROED problem
\eqref{prob:budget-constrained-stochastic-roed}, we require gradient of the
utility function \eqref{eq:average-robust-parameter} with respect to the
uncertain parameter $\robustparam$.

The discussion so far has been agnostic to the dependence of the OED
problem on the uncertain parameter.
However, to differentiate, we will need to explicitly address the uncertain
parameter $\robustparam$ and how it enters the inverse problem.
In general, $\robustparam$ can be a hyperparameter characterizing uncertainty or
misspecification in one or more elements of the inverse problem such as the observation
error model, the prior, or the simulation model.

In our work, as a simplifying assumption, we assume that $\robustparam$ is located
within the observation error covariance $\noisecov$.
That is, $\noisecov \equiv \noisecov(\robustparam)$, where
$\noisecov(\robustparam)$ is assumed to be positive definite for all
$\robustparam \in \robustspace$ and smooth with respect to $\robustparam$.
While this assumption of where $\robustparam$ is located is not strictly
necessary and, in fact, can be relaxed, it simplifies the presentation of the
subsequent methods.

Likewise, at this point, we make the dependence of inverse problem on the design
$\design$ explicit.
In particular, in order to configure active sensors, we construct the modified noise
covariance $\weightednoisecov$ as
\begin{equation} \label{eq:weighted-noise-covariance}
	\weightednoisecov(\design, \robustparam)
	=
	\diag(\design) \noisecov(\robustparam) \diag(\design)
	\, ,
\end{equation}
where $\diag(\design) \in \R^{\Ndata \times \Ndata}$ is a diagonal matrix with the
elements of $\design$ on the diagonal. Likewise, in place of $\noisecovinv$, we use
$\weightednoisecov\pinv(\design, \robustparam)$, where $\dagger$ denotes the
Moore-Penrose pseudoinverse. See \cite{Attia_Constantinescu_2022} for additional details
on how the binary design $\design$ affects the forward model and inverse problem.

Returning to the question of differentiating $\utility$, note that
\begin{equation}
	\grad_{\robustparam} \utility(\design, \robustparam)
	=
	\frac{1}{\Nsaa} \sum_{i=1}^{\Nsaa}
	\grad_{\robustparam} \hat{\utility}(\obs_i, \design, \robustparam) \,,
\end{equation}
it is enough to understand how to differentiate $\hat{\utility}$, the fixed MAP
point approximation to the low-rank information gain.

We next consider the differentiation of $\hat{\utility}$ with respect to
$\robustparam$.
Note, the third term in $\hat{\utility}$ \eqref{eq:low-rank-fixed-map-ig} is
independent of $\robustparam$, hence,
\begin{equation}\label{eqn:utility_grad_robust_param}
	\grad_{\robustparam} \hat{\utility}(\obs_i, \design, \robustparam)
	=
	\grad_{\robustparam} \left(
	\frac{1}{2} \sum_{n=1}^r
	\left[
			\log(1 + \lambda_n^i(\design, \robustparam))
			- \frac{\lambda_n^i(\design, \robustparam)}{1+\lambda_n^i(\design, \robustparam)}
			\right]
	\right) \,.
\end{equation}

An analytical form of the gradient \eqref{eqn:utility_grad_robust_param} can be obtained
by employing an adjoint-based eigenvalue sensitivity framework
\cite{Chowdhary_Tong_Stadler_Alexanderian_2024,Chen_Villa_Ghattas_2019}.
This is done by constructing a Lagrangian over the eigenvalues of systems constraining
the data-misfit Hessian action its eigenproblem.
To perform this technique, we assume that $\ppHm$ is differentiable with respect to
$\robustparam$ and that its dominant eigenvalues are distinct, which is a sufficient
condition for the differentiability of the eigenvalues \cite{Lax_2007}.
For the sake of simplicity of notation, we suppress the data and design dependence of
the Hessian and eigenvalues.

To facilitate the discussion of derivative computation, we consider
\begin{subequations} \label{eq:low-rank-info-gain-constraints}
	\begin{equation}\label{eq:low-rank-info-gain-constraints-objective}
		\frac{1}{2} \sum_{n=1}^r \left[
			\log(1 + \lambda_n) - \frac{\lambda_n}{1+\lambda_n}
			\right] \,,
	\end{equation}
	such that the following eigenproblem constraints hold
	\begin{align}
		\label{eq:eigenproblem-constraint-1}
		\inp{\phi}{\Hm \psi_n}
		 & = \lambda_n \inp{\phi}{\psi_n}_{\priorcovinv},
		 & \forall \phi \in \testspace, \forall n = 1, \ldots, r \,, \\
		\label{eq:eigenproblem-constraint-2}
		\inp{\psi_n}{\psi_n}_{\priorcovinv}
		 & = 1,
		 &                                                           
		\forall n = 1, \ldots, r \,,
	\end{align}
	where $\psi_n = \priorcovsqrt \omega_n$, $\omega_n$ is the eigenvector associated
	with the eigenvalue $\lambda_n$.
	Recall, in \cref{eq:gn-data-misfit-hessian-system}, we stated adjoint-based
	expressions for the action of $\Hm$ in terms of the weak form
	$\weakpde(\state,\invparam,\test)$.
	Restating it for clarity, we therefore also have the constraints
	\begin{equation}
		\Hm(\invparam)(\psi_n, \phi)
		= \inp{\phi}{
			\weakpde_{\invparam\test}(\state, \invparam, \test) \hat{\test}
		} \,,
	\end{equation}
	with state and adjoint constraints
	\begin{align}
		\inp{\tilde{\test}}{\weakpde_\test(\state, \invparam, \test)}
		 & = 0,
		 & \forall \tilde{\test} \in \testspace  \,,  \\
		\inp{\tilde{\state}}{\weakpde_{\state}(\state, \invparam, \test)}
		+ \inp{\tilde{\state}}{\Obs^*
			\weightednoisecov\pinv(\design, \robustparam)
			(\obs - \Obs\state)}
		 & = 0,
		 & \forall \tilde{\state} \in \statespace \,,
	\end{align}
	and incremental state and adjoint constraints for $n = 1, \ldots, r$:
	\begin{align}
		\inp{\tilde{\test}}{\weakpde_{\test\state}(\state, \invparam, \test) \hatstate_n}
		+ \inp{\tilde{\test}}{\weakpde_{\test\invparam}(\state, \invparam, \test) \psi_n}
		 & = 0,
		 & \forall \tilde{\test} \in \testspace \,,   \\
		\inp{\tilde{\state}}{
			\weakpde_{\state\test}(\state, \invparam, \test) \hattest_n
		}
		+ \inp{\tilde{\state}}{
			\Obs^*
			\weightednoisecov\pinv(\design, \robustparam)
			\Obs \hatstate_n
		}
		 & = 0,
		 & \forall \tilde{\state} \in \statespace \,.
	\end{align}
\end{subequations}

To differentiate through~\cref{eq:low-rank-info-gain-constraints-objective}, we first recognize that we can replace
$\lambda_n$ by $\Hm(\invparam)(\psi_n, \psi_n)$ in
\eqref{eq:low-rank-info-gain-constraints-objective}. This eliminates the constraint
\eqref{eq:eigenproblem-constraint-1}. Additionally, to ease the burden of notation,
henceforth we drop the dependence of $\weakpde$ on $(\state, \invparam, \test)$ and
simply write $\weakpde$.
Therefore, a meta-Lagrangian for \eqref{eq:low-rank-info-gain-constraints} is given by
\begin{equation} \label{eq:low-rank-info-gain-langrangian}
	\begin{aligned}
		 & \mc{L}^{\rm IG}\left(
		\state, \invparam, \test,
		\{\psi_n\}_{n=1}^r,
		\{\hatstate_n\}_{n=1}^r,
		\{\hattest_n\}_{n=1}^r,
		\state^*, \test^*,
		\{\lambda_n^*\}_{n=1}^r,
		\{\hatstate_n^*\}_{n=1}^r,
		\{\hattest_n^*\}_{n=1}^r
		; \robustparam
		\right)
		\\
		 & \quad=
		\frac{1}{2} \sum_{n=1}^r \left[
			\log(1
			+ \inp{\psi_n}{
				\weakpde_{\invparam\test}  \hattest_n
			}
			)
			- \frac{
				\inp{\psi_n}{
					\weakpde_{\invparam\test}  \hattest_n
				}
			}{
				1
				+ \inp{\psi_n}{
					\weakpde_{\invparam\test}  \hattest_n
				}
			}
		\right]                                                                                         \\
		 & \qquad + \inp{\test^*}{\weakpde_\test}
		+ \inp{\state^*}{\weakpde_\state}
		+ \inp{\state^*}{\Obs^* \weightednoisecov\pinv (\obs - \Obs\state)}                             \\
		 & \qquad + \sum_{n=1}^r \left[
			\inp{\hattest_n^*}{
				\weakpde_{\test\state} \hatstate_n
				+ \weakpde_{\test\invparam} \psi_n
			}
			+ \inp{\hatstate_n^*}{
				\weakpde_{\state\test} \hattest_n
			}
			+ \inp{\hatstate_n^*}{\Obs^* \weightednoisecov\pinv \Obs \hatstate_n}
		\right]                                                                                         \\
		 & \qquad + \sum_{n=1}^r \lambda_n^* \left[ \inp{\psi_n}{\psi_n}_{\priorcovinv} - 1 \right] \,.
	\end{aligned}
\end{equation}

Subsequently, we proceed to determine the Lagrange multipliers.
By differentiation with respect to $\hattest_n$ in direction $\tildetest$,
and by setting the result to zero we have
\begin{equation}
	\label{eq:lagrange-multiplier-incremental-state}
	\frac{\lambda_n}{2(1+\lambda_n)^2} \inp{\psi_n}{\weakpde_{\invparam\test} \tildetest}
	+ \inp{\hatstate_n^*}{\weakpde_{\state\test} \tildetest}
	= 0, \quad \forall \tildetest \in \testspace \,.
\end{equation}
Reversing the order of differentiation in each of the inner products
shows
this is a rescaled version of the incremental state equation.
In particular, $\hatstate_n^* = \frac{1}{2} \lambda_n(1+\lambda_n)^{-2}
	\hatstate_n$ solves the above, determining the Lagrange multiplier.
Now, by differentiating with respect to $\hatstate_n$ in direction $\tildestate$ and
by setting the result to zero we obtain:
\begin{equation}
	\label{eq:lagrange-multiplier-incremental-adjoint}
	\inp{\hattest_n^*}{ \weakpde_{\test\state} \tildestate }
	+ \inp{\hatstate_n^*}{\Obs^* \weightednoisecov\pinv \Obs \tildestate}
	= 0, \quad \forall \tildestate \in \testspace \,.
\end{equation}
Again, by reversing the order of differentiation in each of the inner products,
we note that this is a rescaled version of the incremental adjoint equation.
Selecting $\hattest_n^* =
	\frac{1}{2} \lambda_n(1+\lambda_n)^{-2} \hattest_n$ solves the above, determining
the multiplier.
Now, differentiating with respect to $\psi_n$ in direction $\phi$ and setting the result
to zero yields
\begin{equation}
	\inp{\hattest_n^*}{\weakpde_{\test\invparam} \phi}
	+ \lambda_n^* \inp{\phi}{\psi_n}_{\priorcovinv}
	= 0, \quad \forall \phi \in \testspace \,.
\end{equation}
The first three terms form a rescaled Hessian action on $\psi_n$ in direction $\phi$.
This, by the definition of eigenfunctions, is precisely equal to $\lambda_n
	\inp{\phi}{\psi_n}_{\priorcovinv}$ times the rescaling factor.
Thus, the selection of $\lambda_n^* = \frac{1}{2} \lambda_n^2(1+\lambda_n)^{-2}$
satisfies the equation.
Now, let's differentiate with respect to $\test$ in direction $\tildetest$ and set the
result to zero to find:
\begin{equation}
	\begin{aligned}
		 & \sum_{n=1}^r \left[
			\frac{\lambda_n}{2(1+\lambda_n)^2}
			\inp{\tildetest}{
				\inp{\psi_n}{
					\weakpde_{\test\invparam\invparam} \psi_n
					+ \weakpde_{\test\invparam\state} \hatstate_n
					+ \weakpde_{\test\invparam\test}  \hattest_n
				}
			}
		\right]                            \\
		 & \quad+
		\sum_{n=1}^r \left[
			\inp{\tildetest}{
				\inp{\hattest_n^*}{
					\weakpde_{\test\test\state} \hatstate_n
					+ \weakpde_{\test\test\invparam} \psi_n
				}
				+ \inp{\hatstate_n^*}{
					\weakpde_{\test\state\test} \hattest_n
					+ \weakpde_{\test\state\state} \hatstate_n
					+ \weakpde_{\test\state\invparam} \psi_n
				}
			}
			\right]
		\\
		 & \qquad\qquad+ \inp{\tildetest}{
			\weakpde_{\test\test}\test^*
			+ \weakpde_{\test\state}\state^*
		}
		= 0\,,\quad \forall \tildetest \in \testspace \,.
	\end{aligned}
\end{equation}
We can further simplify this by noting that all terms involving two derivatives of
$\weakpde$ with respect to $\test$ vanish.
Hence, we have
\begin{equation}
	\label{eq:lagrange-multiplier-state}
	\begin{aligned}
		 & \sum_{n=1}^r \left[
			\inp{\tildetest}{
				\frac{\lambda_n}{2(1+\lambda_n)^2}
				\inp{\psi_n}{
					\weakpde_{\test\invparam\invparam} \psi_n
					+ \weakpde_{\test\invparam\state} \hatstate_n
				}
				+ \inp{\hatstate_n^*}{
					\weakpde_{\test\state\state} \hatstate_n
					+ \weakpde_{\test\state\invparam} \psi_n
				}
			}
		\right]                                                     \\
		 & \qquad+ \inp{\tildetest}{\weakpde_{\test\state}\state^*}
		= 0\,, \quad \forall \tildetest \in \testspace \,.
	\end{aligned}
\end{equation}
As we already know a value for $\hatstate_n^*$, the above is a fully specified equation
for $\state^*$.
Now, let us differentiate with respect to $\state$ in direction $\tildestate$ and set
the result to zero to find:
\begin{equation}\label{eqn:final_Lag}
	\begin{aligned}
		 & \sum_{n=1}^r \left[
			\frac{\lambda_n}{2(1+\lambda_n)^2}
			\inp{\tildestate}{
				\inp{\psi_n}{
					\weakpde_{\state\invparam\test}  \hattest_n
				}
			}
			\right]
		+ \inp{\tildestate}{
			\weakpde  _{\state\test}\test^*
		}
		+ \inp{\tildestate}{\Obs^* \weightednoisecov\pinv \Obs \state^*}
		\\
		 & \qquad+
		\sum_{n=1}^r \left[
			\inp{\tildestate}{
				\inp{\hattest_n^*}{
					\weakpde_{  \state\test\state} \hatstate_n
					+ \weakpde_{\state\test\invparam} \psi_n
				}
				+ \inp{\hatstate_n^*}{
					\weakpde_{  \state\state\test} \hattest_n
				}
			}
			\right]
		= 0\,,\quad \forall \tildetest \in \testspace \,.
	\end{aligned}
\end{equation}

Note that \eqref{eqn:final_Lag} is fully specified, and therefore the associated Lagrange is
entirely determined.
With this, all Lagrange multipliers are specified.
Hence, we can differentiate $\hat{\utility}$ \eqref{eq:low-rank-fixed-map-ig} in
direction $\robustparam$ to obtain $\grad_{\robustparam} \hat{\utility}$ and thus obtain
the analytical form of the gradient of our ROED utility function
\eqref{eqn:roed_utility_function} as:
\begin{subequations}\label{eqn:roed_utility_function_grad_robust_param}
	\begin{align}
		\label{eq:saa-low-rank-fixed-map-eig-grad-robust-param}
		 & \grad_{\robustparam} \utility(\design, \robustparam)
		=
		\frac{1}{\Nsaa} \sum_{i=1}^{\Nsaa}
		\left[
			\sum_{j=1}^{\Nrobust}
			\pdv{}{\theta_j} \hat{\utility}(\obs_i, \design, \robustparam)
			\mbf{e}_j
			\right]
		\, ,
		%
		\\
		\label{eq:saa-low-rank-fixed-map-ig-grad-robust-param}
		 & \begin{aligned}
			   \pdv{}{\theta_i} \hat{\utility}(\obs_i, \design, \robustparam)
			    & =
			   -\inp{\state^*}{
				   \Obs^*
				   \left[
					   \pdv{}{\theta_i} \weightednoisecov\pinv(\design, \robustparam)
					   \right]
				   (\obs_i - \Obs\state)
			   }                       \\
			    & \qquad- \sum_{n=1}^r
			   \inp{\hatstate^*_n}{
				   \Obs^*
				   \left[
					   \pdv{}{\theta_i} \weightednoisecov\pinv(\design, \robustparam)
					   \right]
				   \hatstate_n
			   }
			   \,.
		   \end{aligned} \\
		\label{eq:weighted-noise-covariance-grad-robust-param}
		 &
		\pdv{}{\theta_i} \weightednoisecov\pinv(\design, \robustparam)
		=
		-\weightednoisecov\pinv(\design, \robustparam)
		\left[
			\pdv{}{\theta_i} \noisecov(\robustparam)
			\right]
		\weightednoisecov\pinv(\design, \robustparam) \,.
	\end{align}
\end{subequations}
For more details regarding the derivation of
\cref{eq:weighted-noise-covariance-grad-robust-param}, see \cite[Appendix
	A]{Attia_Leyffer_Munson_2023}.
Note that out of all Lagrange multipliers specified for this expression only $\state^*$
and $\hatstate_n^*$ are necessary to compute.
A detailed algorithmic procedure for calculating
\eqref{eqn:roed_utility_function_grad_robust_param}, in the context of our proposed ROED
approach, is described by \Cref{alg:low-rank-fixed-map-ig-and-eig}.

\subsection{Algorithmic Statement and Computational Considerations}
\label{subsec:computational-complexity}
The developments in the previous sections provide the building blocks of our ROED
framework for nonlinear inverse problems governed by PDEs.
In this section, we provide a complete summary of the steps and the computational
complexity of the proposed algorithm.
\Cref{alg:stochastic-budget-roed} describes our approach for
solving the budget-constrained probabilistic ROED
problem defined by \Cref{prob:budget-constrained-stochastic-roed}.
As described in \Cref{subsec:oed}, the algorithm proceeds by alternating two
steps at each iteration $l$.
First, the conditional Bernoulli model parameter (the policy) $\policy$ is
updated (Step \ref{algstep:Policy_Update})
by using a stochastic optimization procedure where a finite sample of the uncertain
parameter $\finiterobustspace$ is used.
Second, the uncertain parameter $\robustparam$ is updated (Step \ref{algstep:robust-param-opt})
by following a gradient-based optimization approach,
and the optimal solution is used to expand $\finiterobustspace$.

\begin{algorithm}[!htb]
	\caption{Algorithm for solving the ROED problem \eqref{eq:budget-constrained-stochastic-roed}
		\label{alg:stochastic-budget-roed}
	}
	\begin{algorithmic}[1]
		\Require{Initial policy parameter $\policy^{(0)}\in[0, 1]^{\Ndata}$;
			learning rate $\eta_{\policy}$;
			sample size $\Nens$;
			budget $\Nbudget$;
			and an initial sample
			$\finiterobustspace^{(k)}:=\{\robustparam^{(i)}\in\robustspace\,|\,
				i=1, \ldots, k\}$. 
		}
		\Ensure{
			$\optdesign$
		}
		\State Let $l \gets 0$.
		\While{Not Converged} \label{algstep:polyak-outer-loop}

		\Statex \hspace{1.25em} {\color{gray}$\triangleright$ \textit{Outer Optimization: Policy Update}}

		\State \label{algstep:Policy_Update} \(
		\policy^{(l+1)} \gets
		\textsc{policyOpt}(
		\policy^{(l)},\, \Nbudget,\, \eta_{\policy},\, \finiterobustspace^{(k+l)}
		)
		\)

		\State Sample
		$
			S^{(l+1)} \gets
			\{
			\design_i \sim \Probability(\design \, |\, \policy^{(l+1)}, \budget = \Nbudget)
			| i = 1, \dots, \Nens
			\}
		$ 

		\State \(
		\design^{(l+1)}
		\gets \argmax\limits_{\design \in S^{(l+1)}}
		\min\limits_{\robustparam \in \finiterobustspace^{(k+l)}}
		\utility(\design, \robustparam)
		\)

		\Statex {\color{gray}$\triangleright$ \textit{Inner Optimization: Uncertain Parameter Update}}
		\State \label{algstep:robust-param-opt}\(
		\robustparam^{(k+l)}
		\gets
		\argmin\limits_{\robustparam \in \robustspace}
		\utility(\design^{(l+1)}, \robustparam)
		\) \Comment{Requires $\utility,\, \grad_{\robustparam} \utility$; \, use \Cref{alg:low-rank-fixed-map-ig-and-eig}}
		\State \(
		\finiterobustspace^{(k+l+1)}
		\gets \finiterobustspace^{(k+l)} \cup \{\robustparam^{(k+l)}\}
		\)
		\State $l \gets l + 1$
		\EndWhile

		\Statex {\color{gray}$\triangleright$ \textit{Sample Final (Optimal) Policy and Obtain Robust Optimal Design}}
		\State \(
		S_f \gets
		\{
		\design_i \sim \Probability(\design \,|\, \policy^{(l)}, \budget = \Nbudget)
		| i = 1, \dots, \Nens
		\}
		\) 
		\State
		$
			\optdesign
			\gets
			\argmax\limits_{\design \in S_f}
			\min\limits_{\robustparam \in \finiterobustspace^{(k+l)}}
			\utility(\design, \robustparam)
		$
		\State \Return
		$
			\optdesign 
		$
		\Statex

		\Function{policyOpt}{$\policy^{(0)},\, \Nbudget,\, \eta,\, \robustspace$}
		\State Let $n \gets 0$.
		\While{Not Converged}\label{algstep:policy-opt}
		\State Sample
		$
			S \gets
			\{
			\design_i \sim \Probability(\design \, |\,  \policy, \budget=\Nbudget)
			| i = 1, \dots, \Nens
			\}
		$ \Comment{\cite[Algorithm 3.1]{Attia_2024}}
		\State \label{algstep:optimal_baseline_calculation}
		Compute optimal baseline $\optbaseline$ \Comment{\eqref{eq:stochastic-budget-roed-optimal-baseline}}
		\State \label{algstep:stochastic_gradient_calculation}
		Compute the stochastic gradient
		$\grad_{\policy} \stochobj^{\optbaseline}(\policy^{(n)})$ \Comment{\eqref{eq:stochastic-roed-gradient-baseline}}

		\Statex {\color{gray}$\triangleright$ \textit{$\proj$ is a box-constraint projector onto $[0, 1]^{\Ndata}$; see e.g., \cite[Section 4.1]{Attia_2024}}}
		\State \(
		\policy^{(n+1)}
		\gets \policy^{(n)}
		+ \eta P\left(\grad_{\policy} \stochobj^{\optbaseline}\left(\policy^{(n)}\right)\right)
		\) 
		\State $n \gets n + 1$
		\EndWhile
		\State \Return $\policy^{(n)}$
		\EndFunction
	\end{algorithmic}
\end{algorithm}

Note that \Cref{alg:stochastic-budget-roed} is valid for any choice of the utility function
$\utility$ and the method used to evaluate and differentiate it with respect to the
uncertain parameter $\robustparam$.
In our work, we've developed the machinery needed to evaluate and differentiate the
low-rank EIG with a fixed MAP point approximation \eqref{eqn:roed_utility_function}.
Here, we provide \Cref{alg:low-rank-fixed-map-ig-and-eig} to fully specify these techniques.

\begin{algorithm}
	\caption{
		Low-Rank EIG with Fixed MAP Approximation
	}
	\label{alg:low-rank-fixed-map-ig-and-eig}
	\begin{algorithmic}[1]
		\State Compute
		$\avgrobustparam=\frac{1}{|\finiterobustspace|}\sum_{\robustparam \in \finiterobustspace} \robustparam$
		\Comment{\eqref{eq:average-robust-parameter}}

		\State Sample $\{\invparam_i \sim \prior \, |\,\, i = 1, \ldots, \Nsaa\}$ and $\{\noise_i
			\sim \mc{N}(\vzero,\, \noisecov(\avgrobustparam)) \,|\,\, i = 1, \ldots, \Nsaa\}$.

		\State Compute \(
		\{\obs_i = \POM(\invparam_i) + \noise_i \,|\,\, i = 1, \ldots, \Nsaa\}
		\)

		\State Compute \(
		\{
		\postmean^i = \postmean^i(\alldesign, \avgrobustparam)
		\,|\,\, i = 1, \dots, \Nsaa
		\}
		\)
		\Comment{\eqref{eq:postmean}}

		\Statex{}

		\Function{$\hat{\utility}$}{$\obs_i, \design, \robustparam$}
		\State Compute the dominant eigenvalues $\{\lambda_n^i\}_{n=1}^r$ of
		$\ppHm^i(\design, \robustparam)$.
		\State \Return \(
		\frac{1}{2} \sum_{n=1}^r
		\left[
			\log(1 + \lambda_n^i)
			- \frac{\lambda_n^i}{1+\lambda_n^i}
			\right]
		+ \frac{1}{2}\norm{\postmean^i - \priormean}_{\priorcovinv}^2
		\)
		\EndFunction
		\Statex{}
		\Function{$\grad_{\robustparam}\hat{\utility}$}{$\obs_i, \design, \robustparam$}

		\Statex {
			\textcolor{gray}{
				$\quad\triangleright$
				\textit{
					For all following, evaluate at parameters:
				}
				$\postmean^i,\,  \obs_i, \, \design,\, \robustparam$
			}}

		\State Compute dominant eigenpairs $\{\lambda_n^i, \omega_n^i\}_{n=1}^r$ of
		$\ppHm^i(\design, \robustparam)$

		%

		\State Solve for state $\state$ and adjoint $\test$
		\Comment{\cref{eq:state,eq:adjoint}}

		\State Let $\psi_n = \priorcovsqrt \omega_n$ for $n=1,\ldots,r$.

		\State Solve for incremental state $\{\hatstate_n \,|\,\, n = 1, \dots, r\}$
		\Comment{\eqref{eq:lagrange-multiplier-incremental-state}}

		\State Let $\hatstate_n^* = \frac{1}{2} \lambda_n (1+\lambda_n)^{-2} \hatstate_n$

		\State Solve for $\state^*$
		\Comment{\eqref{eq:lagrange-multiplier-state}}

		\State \Return the gradient
		$\grad_{\robustparam} \hat{\utility}(\obs_i, \design, \robustparam)$
		\Comment{\eqref{eq:saa-low-rank-fixed-map-ig-grad-robust-param}}

		\EndFunction
		\Statex{}
		\Function{$\utility$}{$\design, \robustparam$}
		\State \Return \(
		\frac{1}{\Nsaa}
		\sum_{i=1}^{\Nsaa}
		\hat{\utility}(\obs_i, \design, \robustparam)
		\)
		\EndFunction
		\Statex{}
		\Function{$\grad_{\robustparam}\utility$}{$\design, \robustparam$}
		\State \Return \(
		\frac{1}{\Nsaa}
		\sum_{i=1}^{\Nsaa}
		\grad_{\robustparam} \hat{\utility}(\obs_i, \design, \robustparam)
		\)
		\EndFunction
	\end{algorithmic}
\end{algorithm}

In the rest of \Cref{subsec:computational-complexity}, we provide a high-level
discussion of the computational cost of the proposed approach.
Specifically, in \Cref{subsubsec:budget-constrained-stochastic-roed-complexity} we discuss
the  overall complexity of \Cref{alg:stochastic-budget-roed} in terms of the number of
evaluations of the utility function $\utility$.
Then, in \cref{subsubsec:utility-complexity} we summarize the computational cost of
evaluating the utility function $\utility$ and its gradient
$\nabla_{\robustparam}\utility$, required for solving the inner optimization problem,
in terms of the number of PDE solves.

\subsubsection{Overall Complexity of Probabilistic ROED }
\label{subsubsec:budget-constrained-stochastic-roed-complexity}
\Cref{alg:stochastic-budget-roed} is conceptual, and one has to specify a stopping
criterion in Step \ref{algstep:polyak-outer-loop}.  A general approach is to use
a combination of maximum number of iterations and/or projected gradient
tolerance.
For both simplicity and clarity,
we discuss the number of utility function evaluations at each iteration
$l$ of \Cref{alg:stochastic-budget-roed}.
Specifically, we discuss the cost of each of the two alternating steps, namely,
the outer (Step \ref{algstep:Policy_Update}) and the
inner (Step \ref{algstep:robust-param-opt}) optimization steps.

%
%

\boldparagraph{The outer optimization: policy update}
The policy update (Step \ref{algstep:Policy_Update} of \Cref{alg:stochastic-budget-roed})
requires evaluating the stochastic gradient
\eqref{eq:stochastic-roed-gradient-baseline} and the associated
optimal baseline estimate \eqref{eq:stochastic-budget-roed-optimal-baseline}
where $\robustparam^{*}[i]$ is found by solving
\eqref{eq:stochastic-budget-roed-optimal-baseline-robustparam} by enumeration over the
finite sample $\finiterobustspace^{(l)}$.
Thus, at iteration $l$ of \Cref{alg:stochastic-budget-roed},
the optimal baseline (Step \ref{algstep:optimal_baseline_calculation})
requires $\Nens \cdot |\finiterobustspace^{(l)}|$
evaluations of $\utility$, where $|\finiterobustspace^{(l)}|$ is the cardinality
of the finite sample $\finiterobustspace^{(l)}$.
The stochastic gradient (Step \ref{algstep:stochastic_gradient_calculation}) reuses
the same values of  $\utility$ evaluated over the finite sample $\finiterobustspace^{(l)}$,
and thus does not require additional function evaluations.
Moreover, the projection of the stochastic gradient does not require evaluations of $\utility$.
Hence, the total number of utility function evaluations required by the outer optimization
(Step \ref{algstep:polyak-outer-loop})
at iteration $l$ is $\Nens \cdot |\finiterobustspace^{(l)}|$.
The cost of each function evaluation is addressed in \Cref{subsubsec:utility-complexity}.

Note as the algorithm iterates, the policy begins to converge and often the same designs
will be sampled multiple times.
Exploiting this can drastically reduce the number of utility evaluations required by
caching the results of previous evaluations.
This is particularly relevant when the policy degenerates such that its entries are
close to zero or one.
This behavior is often observed in practice for problems with unique optimal solutions.
Thus, the computational cost stated above is in fact an upper bound that tends to reduce as the
optimization algorithm proceeds as discussed in the numerical experiments discussed
in \Cref{sec:numerical-experiments}.

\boldparagraph{The inner optimization: uncertain parameter update}
\Cref{alg:stochastic-budget-roed} solves the max-min optimization problem over an
expanding finite sample of the uncertain parameter $\finiterobustspace^{(l)}$.
At iteration $l$ of the algorithm, the inner optimization (Step \ref{algstep:robust-param-opt})
updates the finite sample by adding the solution of the inner minimization problem
over the continuous space $\robustspace$.
Thus, the number of evaluations of utility function and its gradient depend on the
numerical optimization method used to minimize $\utility$.
Given that the utility function and its gradient both require expensive PDE simulations,
we focus in the rest of this section on the computational cost of evaluating the utility
function $\utility$ and its gradient $\grad_{\robustparam}\utility$ in terms of number
of PDE solves.

\subsubsection{Complexity of Evaluating and Differentiating the Utility}
\label{subsubsec:utility-complexity}
%
The number of PDE solves required for evaluating the utility function and its gradient
has a major impact on the overall cost of the proposed approach.
These solves 
can vary significantly depending on factors such as the discretization method, equation
characteristics, and chosen discretization scheme.
Hence, we provide a summary of the computational complexity in terms of the number of
PDE solves required.
The computational cost is summarized by \Cref{tab:computational-complexity-pde-solves}
and is discussed next.

\begin{table}[!htb]
	\centering
	\caption{Computational complexity in the number of PDE solves.}
	\label{tab:computational-complexity-pde-solves}
	\begin{tabular}{ll}
		\toprule
		Procedure                   & Cost (in PDE solves)                \\
		\midrule
		Evaluation                  & $\mc{O}(4\Nbudget \cdot \Nsaa)$     \\
		Gradient                    & $\mc{O}((3+5\Nbudget) \cdot \Nsaa)$ \\
		Simultaneous Value/Gradient & $\mc{O}((3+5\Nbudget) \cdot \Nsaa)$ \\
		\bottomrule
	\end{tabular}
\end{table}

\boldparagraph{Cost of evaluating $\utility$}
First, let us consider the evaluation of $\utility$ \eqref{eqn:roed_utility_function}.
As the computation of the fixed MAP points is performed once offline, we omit them
from the tabulation.
Thus, for each data sample $\{\obs_i\}_{i=1}^{\Nsaa}$, the only work required is an
eigendecomposition of the corresponding prior preconditioned data misfit Hessian
$\ppHm$.
While, in practice, the rank of $\ppHm$ can vary across $\{\obs_i\}_{i=1}^{\Nsaa}$, we
follow a conservative approach and assume that it has the maximum possible rank.
As described in \Cref{subsubsec:low-rank-and-fixed-map-approximation}, during
\cref{alg:stochastic-budget-roed} this is $\Nbudget$.
A standard randomized eigendecomposition technique, for example
see \cite{Halko_Martinsson_Tropp_2011}, therefore will require $\mc{O}(\Nbudget)$
evaluations of $\ppHm$.
As a single evaluation of $\ppHm$ requires solving the state \cref{eq:state}, adjoint
\cref{eq:adjoint}, and both incremental state \cref{eq:gn-incr-state} and adjoint
\cref{eq:gn-incr-adjoint} equations,
the total number of PDE solves required for an evaluation of the utility function
$\utility$ \eqref{eqn:roed_utility_function} is $\mc{O}(4\Nbudget \cdot \Nsaa)$.

As an aside, we note that the state equation is typically the most expensive PDE to
solve in the evaluation of $\utility$.
Indeed, in the context of this work, the state equation is a nonlinear PDE.
Leveraging, for example, a finite element framework, after discretization one may need
to employ an expensive iterative solver for the resulting nonlinear system of equations.
However, the adjoint, incremental state, and incremental adjoint equations are
linear and, therefore, can be solved more efficiently.

\boldparagraph{Cost of evaluating $\grad_{\robustparam}\utility$}
Now, regarding computing $\grad_{\robustparam}\utility$, for each data realization
$\{\obs_i\}_{i=1}^{\Nsaa}$, we again need to low-rank the prior preconditioned data
misfit Hessian.
As before, this incurs a cost of $\mc{O}(4\Nbudget)$.
Then, for each dominant eigenfunction we need to solve the incremental state equation,
requiring an additional $\mc{O}(\Nbudget)$ PDE solves.
Finally, we include the $3$ PDE solves necessary for $\state, \test$, and $\state^*$.
Hence, the cost of evaluating the gradient $\grad_{\robustparam}\utility$
\eqref{eqn:roed_utility_function_grad_robust_param}
is $\mc{O}((3 + 5\Nbudget) \cdot \Nsaa)$ PDE solves.

Note that if both the evaluation of the utility function and its gradient
are required at the same time, the eigendecomposition of the Hessian can
be shared between the two procedures.
In total, that optimization would save $\mc{O}(4\Nbudget\cdot\Nsaa)$ PDE solves from
the cost of doing both operations separately.
As the eigendecomposition is the primary expense in both operations, this can be a
significant savings.
An example of a case where this would be beneficial is in optimization techniques that
operate off of value and gradient pairs, such as L-BFGS.

\section{Numerical Experiments}
\label{sec:numerical-experiments}
We elaborate our approach for an inverse problem constrained by the following elliptic
PDE model
\begin{equation}\label{eq:poisson}
	\begin{aligned}
		-\nabla \cdot (\exp(m) \nabla u)
		 & =
		0 \qquad \text{in } \Omega := (0,1)^2 \,,                \\
		\exp(m) \nabla u \cdot \vec{n}
		 & =
		0 \qquad \text{on } \Gamma_N := \{0,1\} \times (0,1) \,, \\
		u
		 & =
		g \qquad \text{on } \Gamma_{D} := (0,1) \times \{0,1\}  \,.
	\end{aligned}
\end{equation}
Here, $\Omega = (0,1)^2$, $\Gamma_N$ is the union of the left and right edges, and
$\Gamma_D$ is the union of top and bottom edges;
we let
$g \equiv 0$ on
$(0,1) \times \{0\}$ and $g \equiv 1$ on $(0,1) \times \{1\}$.
%

%
%
%
%
%

Here, we consider the inverse problem of estimating the parameter $m$ in 
\eqref{eq:poisson} from noisy observations of the state variable $u$.
We assume a Gaussian prior $m\sim \mathcal{N}(m_{\rm prior}, \priorcov)$ with mean
$m_{\rm prior} \equiv 0$.
The prior covariance is $\priorcov = \mathcal{A}^{-2}$,
where $\mathcal{A}$ is a differential operator given by the elliptic PDE
\begin{equation}
	\mathcal{A} \invparam
	=
	\begin{cases}
		-\gamma \nabla \cdot (\mat{K} \nabla m) + \delta m & \text{in } \Omega \,,         \\
		\mat{K} \nabla m \cdot \vec{n} + \beta m           & \text{on } \partial\Omega \,.
	\end{cases}
\end{equation}
This describes a commonly used approach for defining the prior in study of
infinite dimensional Bayesian inverse problems
\cite{Stuart_2010,Buithanh_Ghattas_Martin_Stadler_2013}.
The hyperparameters $\gamma$ and $\delta$ are such that $\delta\gamma$ govern
the variance of the samples and $\gamma / \delta$ govern the correlation length.
Here, $\beta = \sqrt{\gamma\delta/2}$ is an empirically selected Robin
coefficient chosen to reduce boundary artifacts, as discussed in
\cite{Daon_Stadler_2018}.
Finally, $\mat{K}$ is a symmetric positive definite matrix.
For our experiment, we select $(\gamma, \delta) = (0.04, 0.2)$ and let $\mat{K}
	= \begin{bmatrix} 1.25 & 0.75 \\ 0.75 & 1.25\end{bmatrix}$.

For the utility function \cref{eqn:roed_utility_function}, we select $\Nsaa =
	32$ and pre-compute the $\Nsaa$ MAP points using prior samples generated from
the above prior.
Additionally, for the optimization procedure, we use the stochastic gradient
method detailed in \cref{subsec:budget-constrained-stochastic-roed} for the
outer stage and the L-BFGS-B method for the inner stage \cite{Liu_Nocedal_1989}.
The convergence criteria for \Cref{alg:stochastic-budget-roed} is determined by
achieving a tolerance $\varepsilon = 10^{-12}$ in its projected gradient
\(
\| P( \grad_\policy \stochobj(\policy^{(n)}) ) \| < \e
\),
or by reaching a maximum number of $100$ iterations.
We start \Cref{alg:stochastic-budget-roed}
with
$\scalarpolicy_i^{(0)} = 0.5$ for $i \in \{1, \dots, \Ndata\}$,
initialize $\finiterobustspace^{(0)}$ with $64$ uniformly randomly sampled
points from $\robustspace$, and a default step size of $\eta = 0.25$.

The numerical discretization of the required weak forms is done in
the finite element package \texttt{FEniCS} \cite{Logg_Mardal_Wells_2012}.
We employ linear triangular continuous Galerkin finite elements on a mesh with
4225 degrees of freedom ($64 \times 64$).
The Bayesian inverse problem itself is implemented in the computational inverse
problem framework \texttt{hIPPYlib} \cite{Villa_Petra_Ghattas_2021}, and
\Cref{alg:stochastic-budget-roed} was built into \texttt{PyOED}
\cite{Chowdhary_Ahmed_Attia_2024} and is publicly available.
Finally, the L-BFGS-B implementation is from \texttt{SciPy} \cite{Scipy_2020}
and miscellaneous array computations are performed in \texttt{NumPy}
\cite{Numpy_2020}.

\subsection{Two Sensor Experiment}
\label{subsec:two_sensors_results}
We begin by considering a low-dimensional experiment for the sake of illustration.
Specifically, we consider only two candidate sensor locations in $\Omega$, one
at $(0.5, 0.25)$ and the other at $(0.5, 0.75)$.
Additionally, we prescribe the following noise covariance matrix
\begin{equation}
	\noisecov(\sigma_1, \sigma_2, \rho)
	=
	\begin{bmatrix}
		\sigma_1^2             & \rho \sigma_1 \sigma_2 \\
		\rho \sigma_1 \sigma_2 & \sigma_2^2
	\end{bmatrix} \, .
\end{equation}
Here, the uncertain parameter vector is given by $\robustparam = (\sigma_1, \sigma_2,
	\rho)$, and we let $\robustspace = [0.05, 0.15]^2 \times [0, 0.99]$. With this, we can
now perform the ROED experiment per \Cref{alg:stochastic-budget-roed}, seeking a design
that is robust against $\robustparam$ as per the utility function $\utility$
\eqref{eqn:roed_utility_function}.
In the present simple example, the set of possible designs are
\[
	\design_1 = (0, 0), \quad
	\design_2 = (1, 0), \quad
	\design_3 = (0, 1), \quad
	\design_4 = (1, 1) \, .
\]
\begin{figure}[!htb]
	\begin{subfigure}[c]{0.48\textwidth}
		\includegraphics[width=\textwidth]{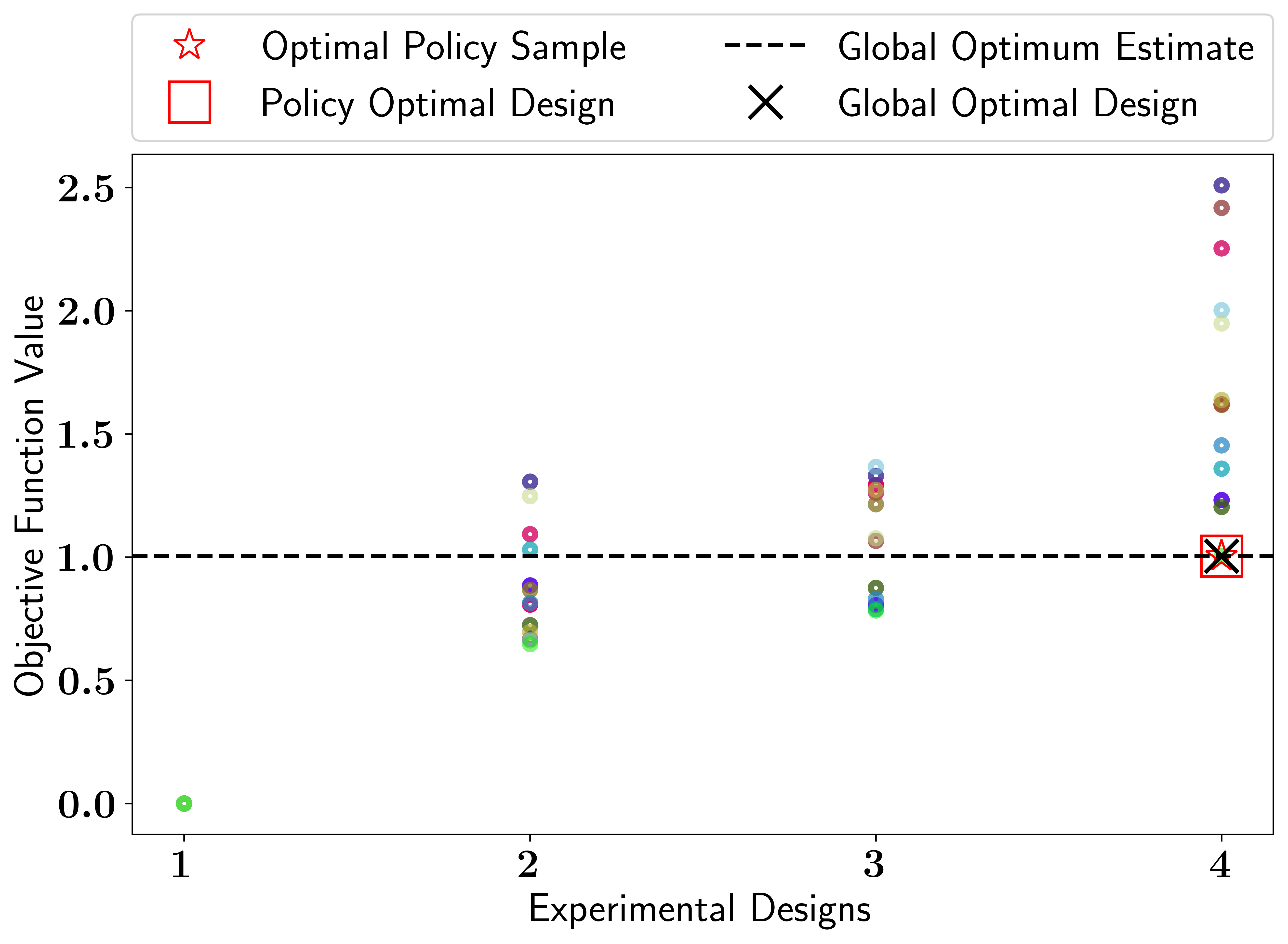}
	\end{subfigure}\hspace{1em}%
	\begin{subfigure}[c]{0.48\textwidth}
		\includegraphics[width=\textwidth]{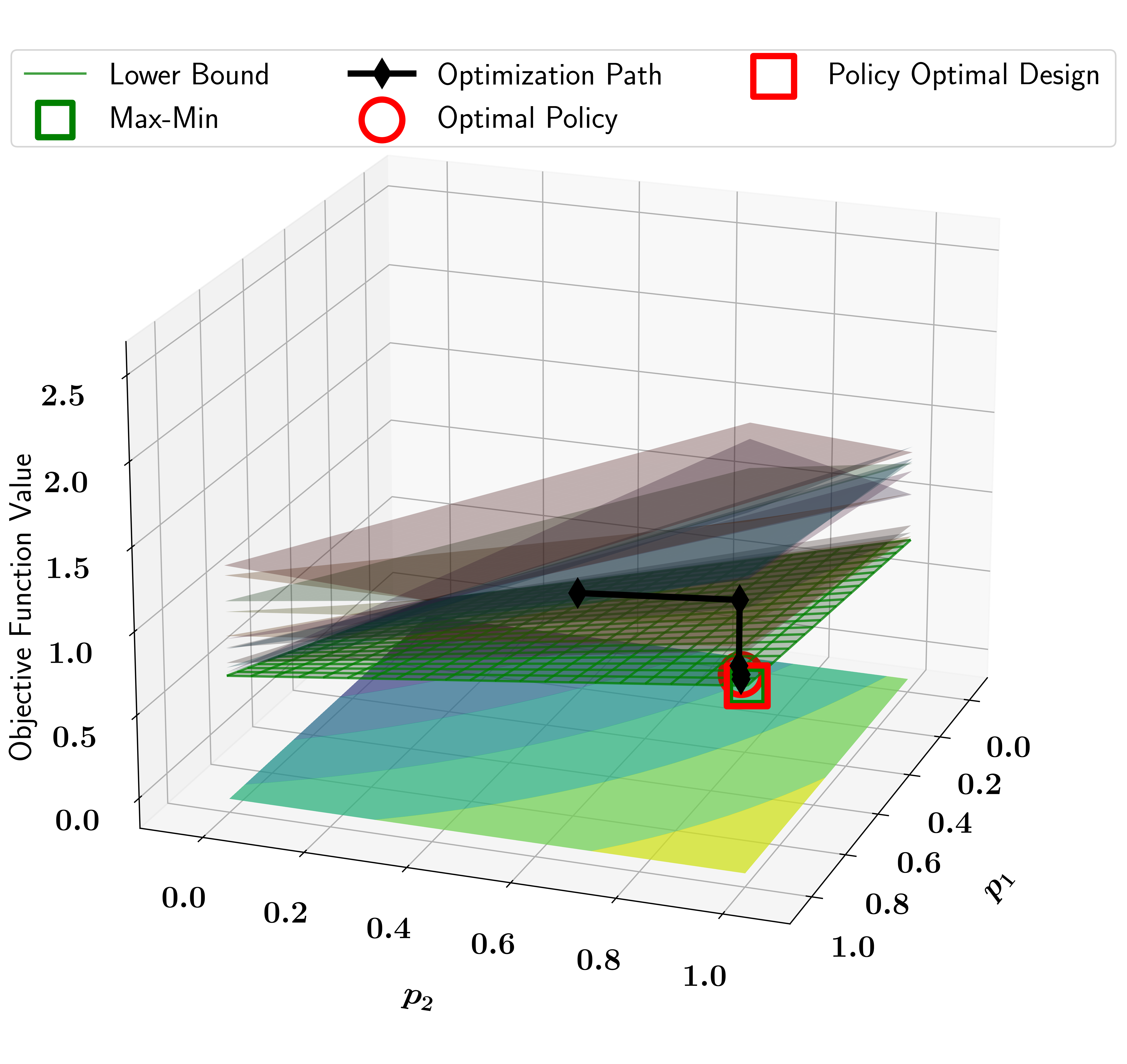}
	\end{subfigure}
	\caption{
		Results of the two sensor experiment. Left: Scatter plot of the utility
		$\utility$ across designs, where different colors represent different
		realizations of $\robustparam$.
		Right: Objective value surfaces where different color surface represents a
		different realization of $\robustparam$.
	}
	\label{fig:two-sensor}
\end{figure}
As we have not imposed a budget, the optimal design is expected to be one with
both sensors active $\design_4$. Indeed, in \Cref{fig:two-sensor}~(left), we see
that the optimal design discovered through \cref{alg:stochastic-budget-roed} is
$\optdesign = \design_4$.
Additionally, we note that the optimal policy $\robustoptpolicy$ degenerates to
the optimal design, as seen in \Cref{fig:two-sensor}~(right).
This is expected in a setting where the global optimal design is unique; see
\cite{Attia_Leyffer_Munson_2023} for more details.

Another interesting feature of this setting is that the utility $\utility$ is
not monotone in the noise parameters. Indeed, as both the scatter plot and the
intersecting objective value surfaces in \Cref{fig:two-sensor} (right) indicate,
different designs have different worst-case $\robustparam$.
This, critically, confirms the necessity for a robust optimization procedure.
After all, if one could determine the worst-case $\robustparam$ a priori, one
could simply select those as a nominal value and perform a standard OED
procedure.

\begin{figure}[!htb]
	\centering
	\includegraphics[width=0.75\textwidth]{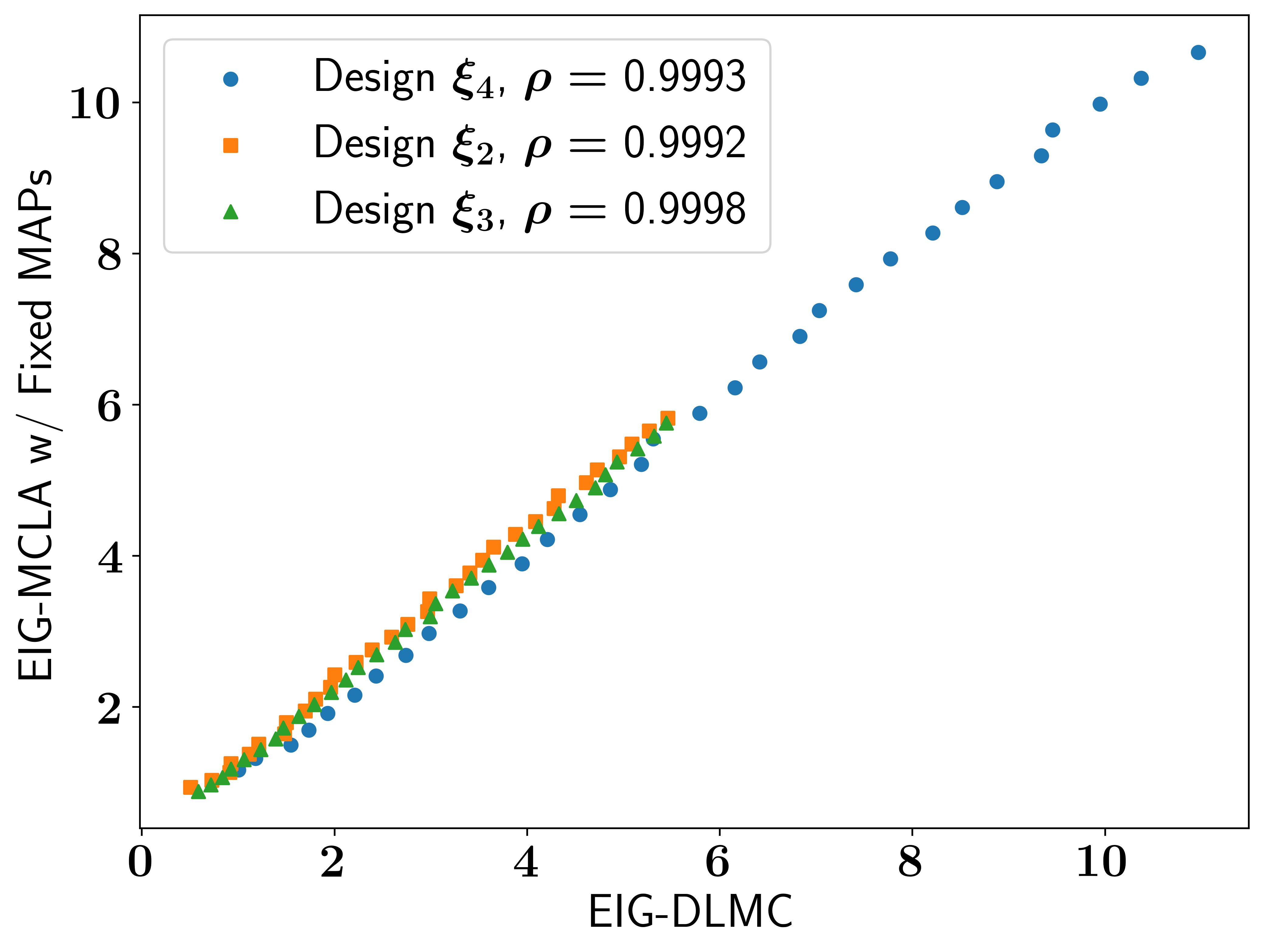}
	\caption{
		Comparison of $\utility$ and a double-loop Monte Carlo estimator of the EIG.
		This is done across different designs and realizations of $\robustparam$.
		Each marker represents a different realization of the noise parameter,
		while the color/shape indicates a different design.
	}
	\label{fig:two-sensor-dlmc-vs-mcla}
\end{figure}

Finally, we assess the approximation quality of $\utility$ by comparing it
to the true EIG as approximated by a double loop Monte Carlo estimator
(\textbf{EIG-DLMC}), see \cite{Huan_Marzouk_2013} for details.
We use $100,000$ samples for the inner integral and $1000$ samples for the outer
integral.
Furthermore, for this experiment we expand $\robustspace$ to include a
wider range of parameters for the sake of a more complete
investigation of both $\utility$ and $\EKLD$ across noise levels.
In particular, we let $\robustspace = [0.01, 0.25]^2 \times [0, 0.99]$.
The results, shown in \Cref{fig:two-sensor-dlmc-vs-mcla}, demonstrate that
across different realizations of the design and noise parameter, both $\utility$
and $\EKLD$ are highly correlated with Pearson correlation coefficients greater
than $0.99$.
Note, we assess the correlation between the two quantities, as opposed to a more
traditional error metric, because we are interested in demonstrating whether the
optimization landscape of both objectives is similar.
The experimental set up for this verification is similar to
\cite{Wu_Chen_Ghattas_2023}, where they considered correlation for a similar
estimator across designs.
Thus, we complement their findings and provide additional insight into the
correlation between $\utility$ and $\EKLD$ for inverse problems of this type.

\subsection{64 Sensor, Budget 8, Experiment}
\label{subsec:64_sensors_results}
Here, we consider an ROED setup with a
realistic number of candidate sensors and a large number of uncertain
parameters.
Specifically, we consider the same simulation model as before, but with $\Ndata=64$
candidate sensor locations forming a regular grid within the domain $\Omega$ and
$\Nbudget = 8$.
Additionally, we consider an observation error covariance matrix
$\noisecov$ whose entries are of the form:
\begin{subequations}
	\begin{equation}
		(\Gamma_\text{n})_{ij}
		=
		\begin{cases}
			\sigma_i^2                                  & \text{if } i = j    \\
			\sigma_i \sigma_j \rho_{ij}(\ell_1, \ell_2) & \text{if } i \neq j
		\end{cases} \, ,
	\end{equation}
	where
	\begin{equation}
		\rho_{ij}(\ell_1, \ell_2)
		= \exp \left(
		-\frac{1}{2\ell_1}| s_i^{(1)} - s_j^{(1)}|
		- \frac{1}{2\ell_2}|s_i^{(2)} - s_j^{(2)}|
		\right) \, ,
	\end{equation}
	with $(s_i^{(1)}, s_i^{(2)})$ denoting the coordinates of the $i$th sensor.
	In this setting, $\robustspace = [0.05, 0.15]^{64} \times [0.01, 2.00]^2$.
\end{subequations}

\begin{figure}[!htb]
	\centering
	\includegraphics[width=\textwidth]{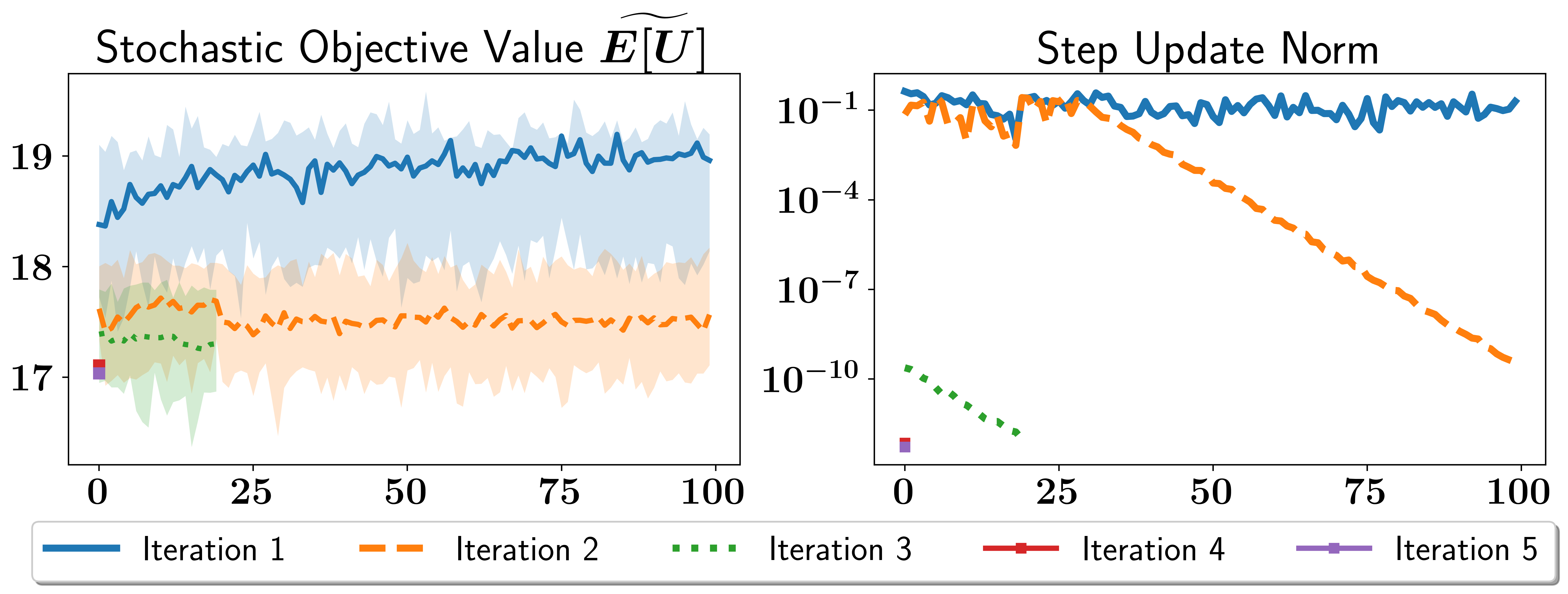}
	\caption{
		Optimization trajectory of the 64 sensor, budget 8, experiment. The iterations refer
		to those of the outer Polyak loop, and the x-axis refers to steps taken in the
		design optimization step of the inner Polyak loop.
		Left: Progress of an estimate to the expectation of the utility $\utility$
		over designs sampled from the policy at that iteration.
		The line represents the mean of the expectation whereas the top and bottom of the
		shaded region represent the maximum and minimum respectively. Right: Norm of the
		update in the policy $\policy$ over the course of the algorithm.
	}
	\label{fig:highdim-roed-progress}
\end{figure}
\begin{figure}[!htb]
	\centering
	\raisebox{-0.6em}{ 
		\includegraphics[width=0.49\textwidth]{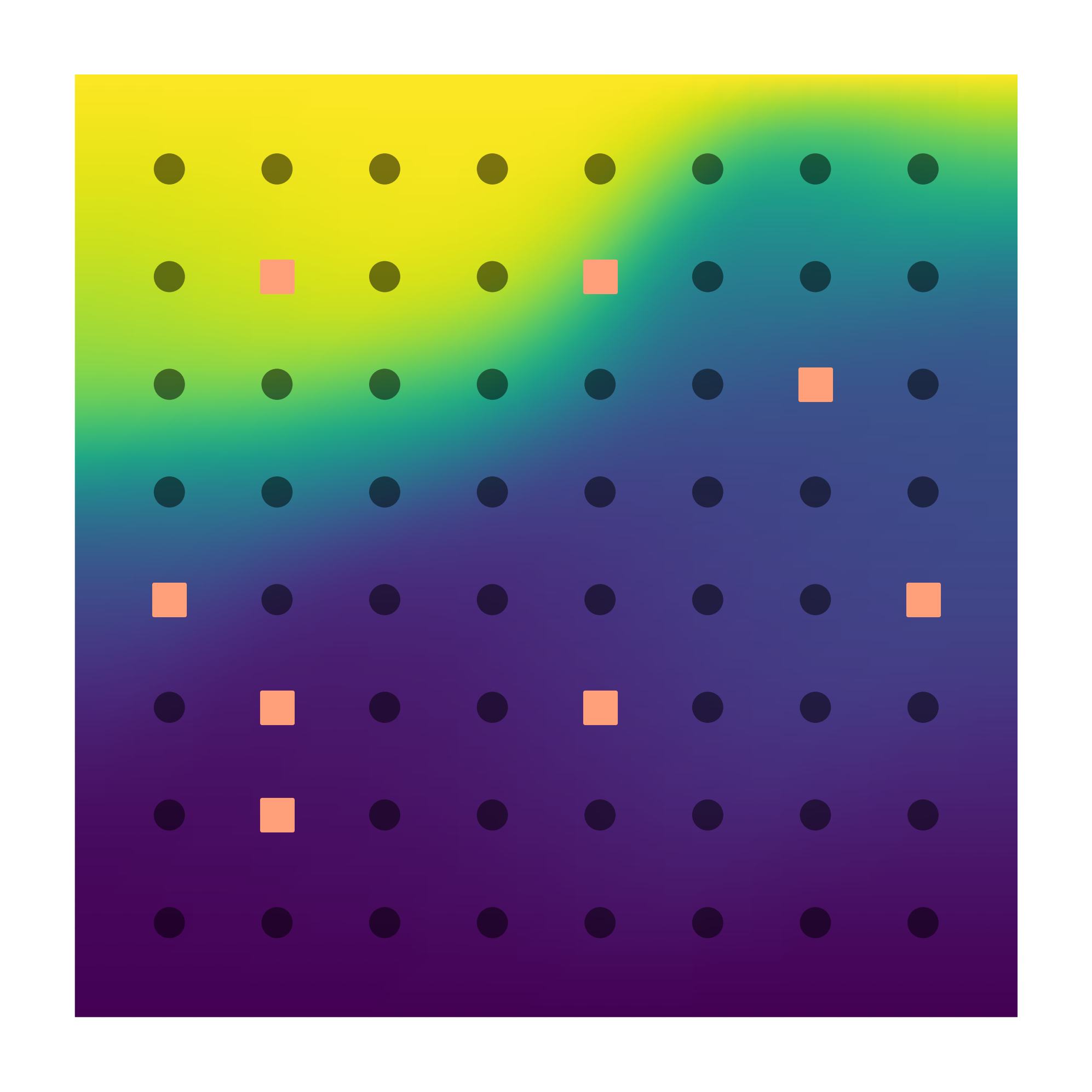}
	}
	\includegraphics[width=0.45\textwidth]{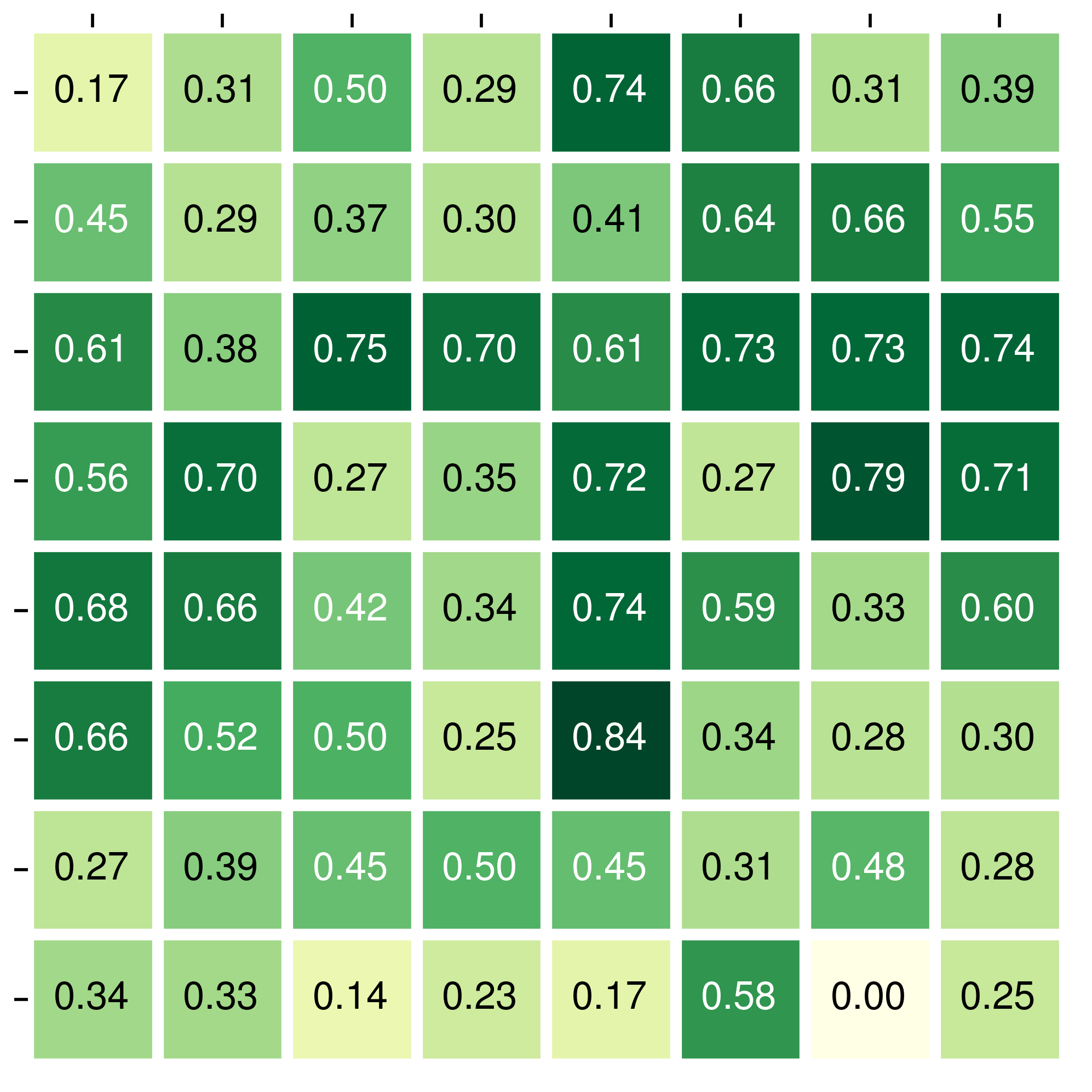}
	\caption{
		Results of the 64 sensor, budget 8, experiment.
		Left: Optimal design discovered by sampling from policy and selecting the design
		with the highest utility.
		Right: Optimal policy $\robustoptpolicy$ discovered by the stochastic optimization
		algorithm, visualized across the sensor grid.
	}
	\label{fig:highdim-optimal-policy-and-design}
\end{figure}
In~\Cref{fig:highdim-roed-progress}, we report the optimization trajectory of
\Cref{alg:stochastic-budget-roed} and in~\Cref{fig:highdim-optimal-policy-and-design} we
show the resulting optimal policy and optimal design.
As demonstrated in \cref{fig:highdim-roed-progress}, the algorithm took five total
iterations of the outer / inner optimization process to converge, though the policy
itself converged early during the third outer optimization stage.
As such, the fourth and fifth outer optimization stages immediately terminated after
a single outer optimization stage iteration.
Note, in this experiment, the optimal policy did not degenerate;
see \Cref{fig:highdim-optimal-policy-and-design} (right).
Likely, this is due to multiple designs producing similar utility.
See the numerical results \cite{Attia_2024} for a similar observation for a different
inverse problem.

\begin{figure}[!htb]
	\begin{subfigure}[c]{0.48\textwidth}
		\includegraphics[width=\textwidth]{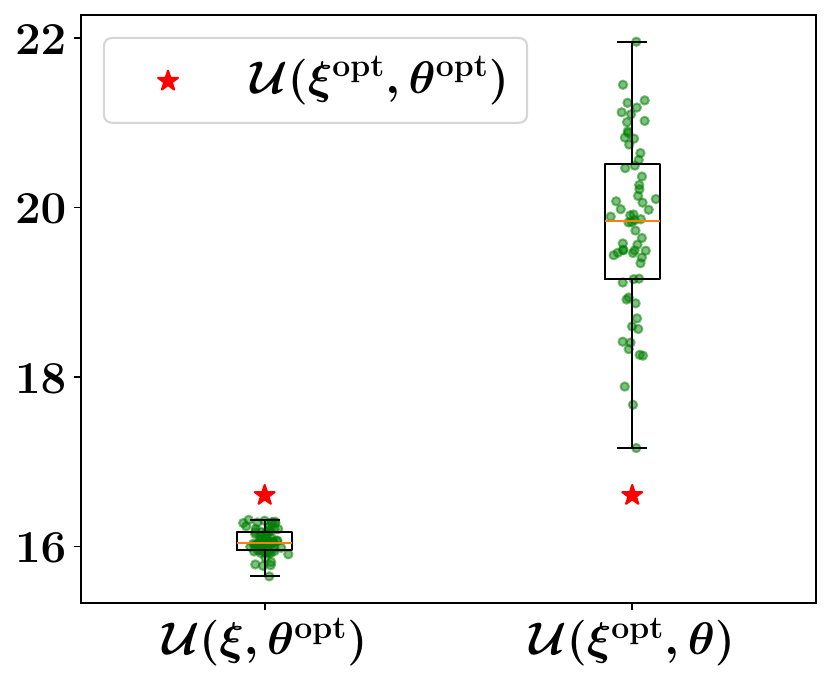}

	\end{subfigure}\hspace{1em}%
	\begin{subfigure}[c]{0.48\textwidth}
		\includegraphics[width=\textwidth]{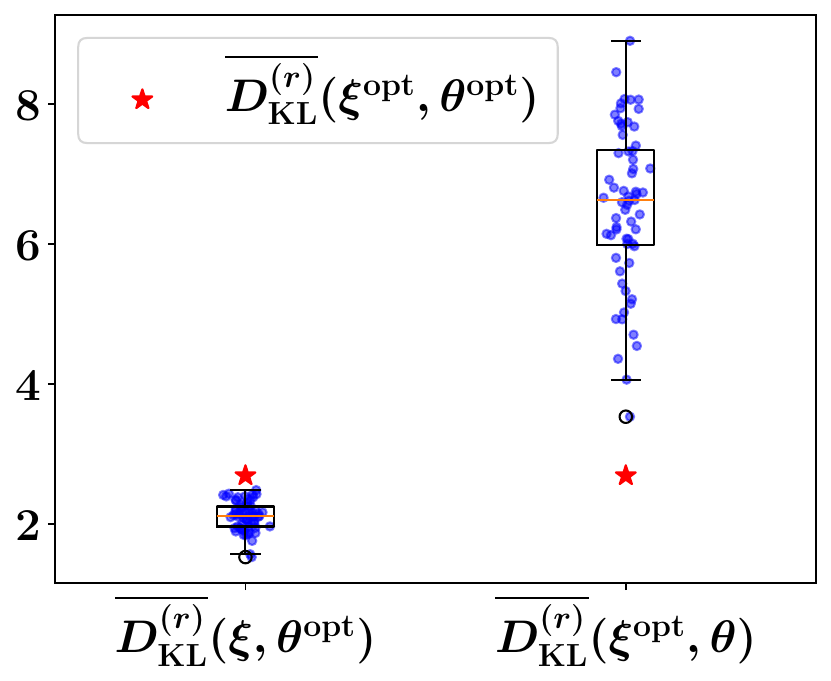}
	\end{subfigure}
	\caption{
		A visualization of the quality of $(\optdesign, \optrobustparam)$ for the 64 sensor,
		budget 8, ROED experiment.
		Here we compare the utility of $(\optdesign, \optrobustparam)$ against $(\optdesign,
			\robustparam)$ for random $\robustparam$ and $(\design, \optrobustparam)$ for random
		$\design$.
		These are evaluated using the utility $\utility$ (left) and the low-rank
		EIG $\lowrankeig$ (right).
	}
	\label{fig:highdim-random-vs-optimal}
\end{figure}
To assess the effectiveness of the proposed strategy, we compare
$\utility(\optdesign, \optrobustparam)$ versus $\utility(\design,
	\optrobustparam)$ for an ensemble of random designs.
While it is not guaranteed that $\utility(\optdesign, \optrobustparam) \geq
	\utility(\design, \optrobustparam)$ for all $\design \in \SNB$, the comparison
is insightful.
We do this comparison in \Cref{fig:highdim-random-vs-optimal}~(left). In the
same figure, we also examine the optimality of the uncertain parameter by
comparing against $\utility(\optdesign, \robustparam)$ for random realizations
of $\robustparam$. For both comparisons, 64 random samples were used.
We see that $\utility(\optdesign, \optrobustparam)$ is significantly higher than
that of $\utility(\design, \optrobustparam)$ for random designs. This indicates
that the discovered design is nearly optimal.
We see similar results for $\optrobustparam$.

In \Cref{fig:highdim-random-vs-optimal}~(right), we repeat the above comparison
using the utility $\lowrankeig$.
The results in the figure indicate that the optimality of $\optdesign$ and
$\optrobustparam$ continue to hold.
This provides numerical evidence regarding the effectiveness of the proposed
strategy and the suitability of the proposed approximation framework.
However, note that the scale of $\utility$ and $\lowrankeig$ are different.

\begin{figure}[!htb]
	\centering
	\includegraphics[width=\textwidth]{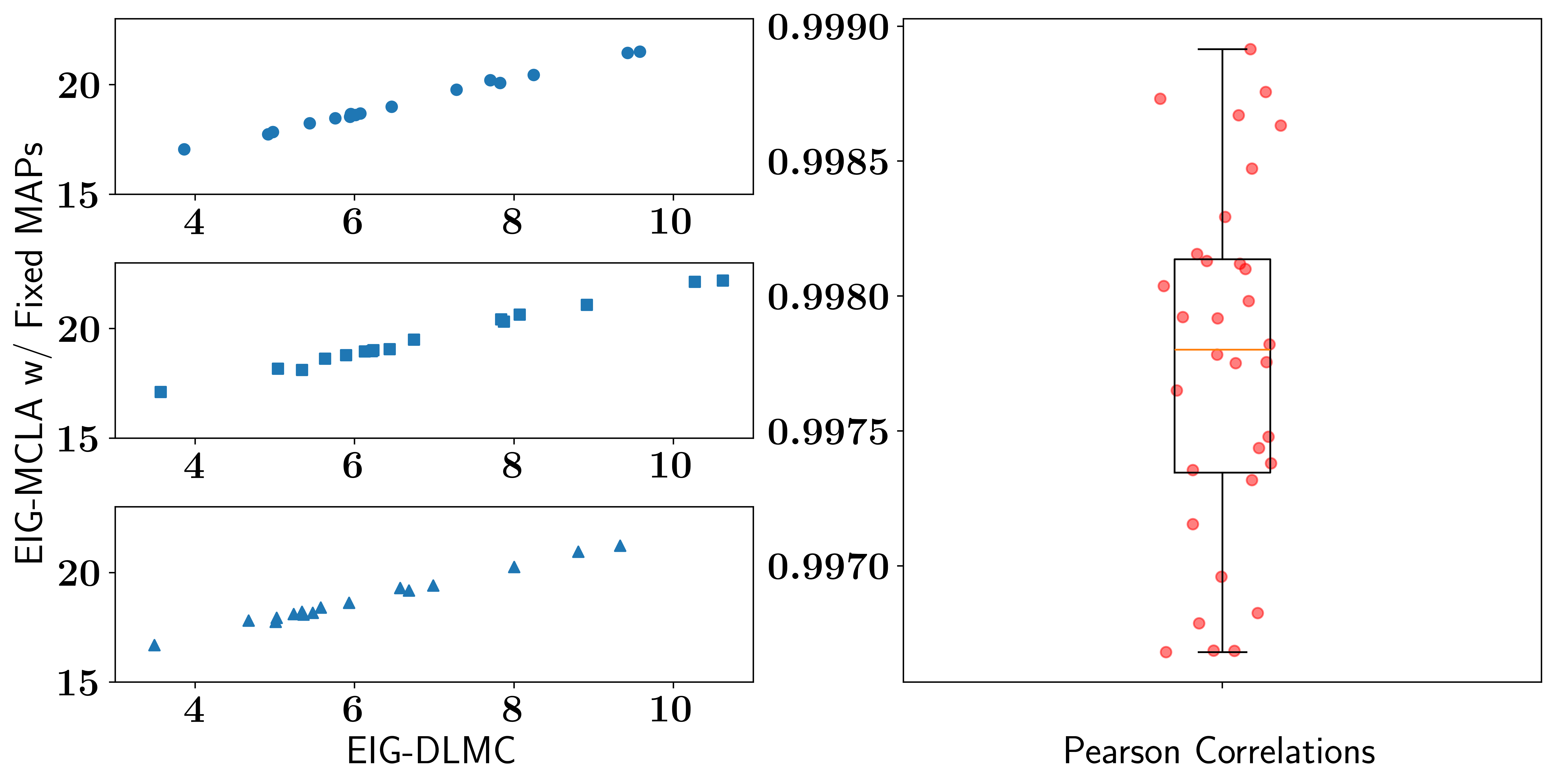}
	\caption{
		Comparison of objective values across different designs and realizations of
		the noise parameter between $\utility$ and an estimate of the true
		expected information gain.
		(Left) Each panel represents a specific randomly selected design, and each
		point represents the utility evaluated at said design across different
		realizations of the noise parameter.
		(Right) A scatter plot of the Pearson correlation coefficients between the
		utility and an approximation to the EIG. Each point represents a fixed
		design and the pearson correlation coefficients are calculated across
		different realizations of the noise parameter.
	}
	\label{fig:highdim-dlmc-vs-mcla}
\end{figure}
Finally, we perform a similar verification as in \Cref{fig:two-sensor-dlmc-vs-mcla} by
assessing the correlation between $\utility$ and the true EIG.
As the dimension of the design space is intractable to enumerate in this
context, we conduct our comparison across $32$ uniformly randomly selected
designs and another $32$ uniformly randomly selected noise parameters.
The results are collected in \Cref{fig:highdim-dlmc-vs-mcla}.
Once again, we consistently observe a strong correlation between the utility and
the true expected information gain.
The figure also illustrates a large bias between the two objectives.
However, as our intention is to produce a surrogate suited for optimization,
a strong correlation positively verifies the similarity of their optimization
landscape.
As mentioned before, this experiment complements those performed in
\cite{Wu_Chen_Ghattas_2023}, and therefore adds to the body of evidence
supporting similar techniques.

\section{Conclusion}
\label{sec:conclusion}

In this article, we have outlined a scalable procedure for robust optimal design
of large-scale Bayesian nonlinear inverse problems
governed by PDEs.
We have constructed a utility agnostic framework for the robust OED problem that
is able to take advantage of the smaller space dictated by the budget.
Then, for the case of the expected information gain used as the utility, we have
developed a framework for its efficient approximation using variational tools in
order to enable its efficient evaluation and differentiation.

There are various avenues for future work. In the first place, we may consider design
criteria other than the expected information gain. Examples include generalizations of
Bayesian A-optimality, goal-oriented design criteria, or decision theoretic ones such as
the Bayes risk.

Secondly, the Laplace approximation and subsequent Gauss-Newton approximation to the
Hessian were critical to the success of the proposed framework. However, while these are
typically appropriate for many problems in practice, there is no guarantee that they
will be appropriate for all problems. Particularly, consider posteriors which are
multimodal or have heavy tails. In this case, an alternative approach may involve
efficient surrogates of the simulation model. There has been significant work in this
area regarding the use of neural networks as surrogates for governing equations given by
PDEs, see \cite{OLearyRoseberry_Villa_Chen_Ghattas_2022}. An investigation into the use
of such surrogates for the robust OED problem is an avenue for future work.

Finally, while the new budget-constrained ROED algorithm eliminates the major drawback
of the penalty parameter tuning stage, there are still several avenues for improvement.
For example, applying a probabilistic approach to the inner optimization stage could
eliminate the need for developing bespoke gradients for the specific selection of
utility and uncertain parameter. Likewise, leveraging different optimization approaches
for outer optimization stage might also prove useful. For example, although relaxed OED
approaches may not be suitable, the use of greedy or exchange type algorithms may prove
to be a simple alternative; see, e.g., \cite{Lau_Zhou_2020}.

\bibliographystyle{siamplain}
\bibliography{references}

\begin{thebibliography}{10}

\bibitem{Alexanderian_2021}
{\sc A.~Alexanderian}, {\em Optimal experimental design for
  infinite-dimensional {Bayesian} inverse problems governed by {PDEs}: a
  review}, Inverse Problems, 37 (2021), p.~043001.

\bibitem{Alexanderian_Gloor_Ghattas_2016}
{\sc A.~Alexanderian, P.~J. Gloor, and O.~Ghattas}, {\em On {Bayesian} {A}- and
  {D}-optimal experimental designs in infinite dimensions}, Bayesian Analysis,
  11 (2016), pp.~671 -- 695.

\bibitem{Alexanderian_Nicholson_Petra_2024}
{\sc A.~Alexanderian, R.~Nicholson, and N.~Petra}, {\em Optimal design of
  large-scale nonlinear {Bayesian} inverse problems under model uncertainty},
  Inverse Problems,  (2024).

\bibitem{Alexanderian_Petra_Stadler_Sunseri_2021}
{\sc A.~Alexanderian, N.~Petra, G.~Stadler, and I.~Sunseri}, {\em Optimal
  design of large-scale {Bayesian} linear inverse problems under reducible
  model uncertainty: Good to know what you don't know}, SIAM/ASA Journal on
  Uncertainty Quantification, 9 (2021), pp.~163--184.

\bibitem{Alexanderian_Saibaba_2018}
{\sc A.~Alexanderian and A.~K. Saibaba}, {\em Efficient {D}-optimal design of
  experiments for infinite-dimensional {Bayesian} linear inverse problems},
  SIAM Journal on Scientific Computing, 40 (2018), pp.~A2956--A2985.

\bibitem{Atkinson_Donev_1992}
{\sc A.~C. Atkinson and A.~N. Donev}, {\em Optimum Experimental Designs},
  Oxford, 1992.

\bibitem{Attia_2024}
{\sc A.~Attia}, {\em Probabilistic approach to black-box binary optimization
  with budget constraints: Application to sensor placement}, arXiv preprint
  arXiv:2406.05830,  (2024).

\bibitem{Attia_Constantinescu_2022}
{\sc A.~Attia and E.~Constantinescu}, {\em Optimal experimental design for
  inverse problems in the presence of observation correlations}, SIAM Journal
  on Scientific Computing, 44 (2022), pp.~A2808--A2842.

\bibitem{Attia_Leyffer_Munson_2023}
{\sc A.~Attia, S.~Leyffer, and T.~Munson}, {\em Robust {A}-optimal experimental
  design for {Bayesian} inverse problems},  (2023),
  \url{https://arxiv.org/abs/arXiv:2305.03855}.

\bibitem{Attia_Leyffer_Munson_2022}
{\sc A.~Attia, S.~Leyffer, and T.~S. Munson}, {\em Stochastic learning approach
  for binary optimization: Application to {Bayesian} optimal design of
  experiments}, SIAM Journal on Scientific Computing, 44 (2022),
  pp.~B395--B427.

\bibitem{Bangerth_2008}
{\sc W.~Bangerth}, {\em A framework for the adaptive finite element solution of
  large-scale inverse problems}, SIAM Journal on Scientific Computing, 30
  (2008), pp.~2965--2989.

\bibitem{Bartuska_Espath_Tempone_2022}
{\sc A.~Bartuska, L.~Espath, and R.~Tempone}, {\em Small-noise approximation
  for {B}ayesian optimal experimental design with nuisance uncertainty},
  Comput. Methods Appl. Mech. Engrg., 399 (2022), p.~115320.

\bibitem{Biedermann_Dette_2003}
{\sc S.~Biedermann and H.~Dette}, {\em A note on maximin and {Bayesian}
  {D}-optimal designs in weighted polynomial regression}, Mathematical Methods
  of Statistics, 12 (2003), p.~358–370.

\bibitem{Buithanh_Ghattas_Martin_Stadler_2013}
{\sc T.~Bui-Thanh, O.~Ghattas, J.~Martin, and G.~Stadler}, {\em A computational
  framework for infinite-dimensional {Bayesian} inverse problems part i: The
  linearized case, with application to global seismic inversion}, SIAM Journal
  on Scientific Computing, 35 (2013), pp.~A2494--A2523.

\bibitem{Chaloner_Verdinelli_1995}
{\sc K.~Chaloner and I.~Verdinelli}, {\em Bayesian experimental design: A
  review}, Statistical Science, 10 (1995).

\bibitem{Chen_Villa_Ghattas_2019}
{\sc P.~Chen, U.~Villa, and O.~Ghattas}, {\em Taylor approximation and variance
  reduction for pde-constrained optimal control under uncertainty}, Journal of
  Computational Physics, 385 (2019), p.~163–186.

\bibitem{Chowdhary_2025}
{\sc A.~Chowdhary}, {\em Scalable Uncertainty Quantification for
  Infinite-Dimensional Bayesian Inverse Problems.}, {PhD} thesis, Feb. 2025,
  \url{https://www.lib.ncsu.edu/resolver/1840.20/45106}.

\bibitem{Chowdhary_Ahmed_Attia_2024}
{\sc A.~Chowdhary, S.~E. Ahmed, and A.~Attia}, {\em {PyOED}: An extensible
  suite for data assimilation and model-constrained optimal design of
  experiments}, ACM Trans. Math. Softw., 50 (2024).

\bibitem{Chowdhary_Tong_Stadler_Alexanderian_2024}
{\sc A.~Chowdhary, S.~Tong, G.~Stadler, and A.~Alexanderian}, {\em Sensitivity
  analysis of the information gain in infinite-dimensional bayesian linear
  inverse problems}, International Journal for Uncertainty Quantification, 14
  (2024).

\bibitem{Cox_1992}
{\sc D.~R. Cox}, {\em Planning of Experiments}, Wiley Classics Library, John
  Wiley \& Sons, Nashville, TN, Apr. 1992.

\bibitem{Daon_Stadler_2018}
{\sc Y.~Daon and G.~Stadler}, {\em Mitigating the influence of the boundary on
  {PDE}-based covariance operators}, 2018.

\bibitem{Darges_Alexanderian_Gremaud_2023}
{\sc J.~Darges, A.~Alexanderian, and P.~Gremaud}, {\em Variance-based
  sensitivity of bayesian inverse problems to the prior distribution},
  International Journal for Uncertainty Quantification.

\bibitem{Dette_Melas_Pepelyshev_2003}
{\sc H.~Dette, V.~B. Melas, and A.~Pepelyshev}, {\em Standardized maximin
  e-optimal designs for the michaelis-menten model}, Statistica Sinica, 13
  (2003), p.~1147–1163.

\bibitem{Federov_2010}
{\sc V.~Fedorov}, {\em Optimal experimental design}, WIREs Computational
  Statistics, 2 (2010), pp.~581--589.

\bibitem{Go_Isaac_2022}
{\sc J.~Go and T.~Isaac}, {\em Robust expected information gain for optimal
  {Bayesian} experimental design using ambiguity sets}, in Proceedings of the
  Thirty-Eighth Conference on Uncertainty in Artificial Intelligence,
  J.~Cussens and K.~Zhang, eds., vol.~180 of Proceedings of Machine Learning
  Research, PMLR, 01--05 Aug 2022, pp.~728--737.

\bibitem{Halko_Martinsson_Tropp_2011}
{\sc N.~Halko, P.~G. Martinsson, and J.~A. Tropp}, {\em Finding structure with
  randomness: Probabilistic algorithms for constructing approximate matrix
  decompositions}, SIAM Review, 53 (2011), pp.~217--288.

\bibitem{Numpy_2020}
{\sc C.~R. Harris, K.~J. Millman, S.~J. van~der Walt, R.~Gommers, P.~Virtanen,
  D.~Cournapeau, E.~Wieser, J.~Taylor, S.~Berg, N.~J. Smith, R.~Kern, M.~Picus,
  S.~Hoyer, M.~H. van Kerkwijk, M.~Brett, A.~Haldane, J.~F. del R{\'{i}}o,
  M.~Wiebe, P.~Peterson, P.~G{\'{e}}rard-Marchant, K.~Sheppard, T.~Reddy,
  W.~Weckesser, H.~Abbasi, C.~Gohlke, and T.~E. Oliphant}, {\em Array
  programming with {NumPy}}, Nature, 585 (2020), pp.~357--362,
  \url{https://doi.org/10.1038/s41586-020-2649-2},
  \url{https://doi.org/10.1038/s41586-020-2649-2}.

\bibitem{Huan_Marzouk_2013}
{\sc X.~Huan and Y.~M. Marzouk}, {\em Simulation-based optimal bayesian
  experimental design for nonlinear systems}, Journal of Computational Physics,
  232 (2013), p.~288–317, \url{https://doi.org/10.1016/j.jcp.2012.08.013}.

\bibitem{Kaipio_Kolehmainen_13}
{\sc J.~Kaipio and V.~Kolehmainen}, {\em Approximate marginalization over
  modeling errors and uncertainties in inverse problems}, Bayesian Theory and
  Applications,  (2013), pp.~644--672.

\bibitem{Kolehmainen_Tarvainen_Arridge_EtAl_11}
{\sc V.~Kolehmainen, T.~Tarvainen, S.~R. Arridge, and J.~P. Kaipio}, {\em
  Marginalization of uninteresting distributed parameters in inverse
  problems-application to diffuse optical tomography}, International Journal
  for Uncertainty Quantification, 1 (2011).

\bibitem{Koval_Alexanderian_Stadler_2020}
{\sc K.~Koval, A.~Alexanderian, and G.~Stadler}, {\em Optimal experimental
  design under irreducible uncertainty for linear inverse problems governed by
  {PDEs}}, Inverse Problems, 36 (2020).

\bibitem{Kullback_Leibler_1951}
{\sc S.~Kullback and R.~A. Leibler}, {\em On information and sufficiency}, The
  Annals of Mathematical Statistics, 22 (1951), pp.~79--86.

\bibitem{Lau_Zhou_2020}
{\sc L.~C. Lau and H.~Zhou}, {\em A local search framework for experimental
  design},  (2020).

\bibitem{Lax_2007}
{\sc P.~D. Lax}, {\em Linear algebra and its applications}, Pure and Applied
  Mathematics: A Wiley Series of Texts, Monographs and Tracts, Wiley-Blackwell,
  Chichester, England, 2~ed., Aug. 2007.

\bibitem{Levitin_Polyak_1966}
{\sc E.~Levitin and B.~Polyak}, {\em Constrained minimization methods}, USSR
  Computational Mathematics and Mathematical Physics, 6 (1966), pp.~1--50.

\bibitem{Liu_Nocedal_1989}
{\sc D.~C. Liu and J.~Nocedal}, {\em On the limited memory {BFGS} method for
  large scale optimization}, Mathematical Programming, 45 (1989), p.~503–528.

\bibitem{Logg_Mardal_Wells_2012}
{\sc A.~Logg, K.-A. Mardal, and G.~Wells}, {\em Automated Solution of
  Differential Equations by the Finite Element Method: The {FEniCS} Book},
  vol.~84 of Lecture Notes in Computational Science and Engineering, Springer
  Berlin Heidelberg, Berlin, Heidelberg, 2012.

\bibitem{Mozumder_Tarvainen_Arridge_EtAl_16}
{\sc M.~Mozumder, T.~Tarvainen, S.~Arridge, J.~P. Kaipio, C.~D'Andrea, and
  V.~Kolehmainen}, {\em Approximate marginalization of absorption and
  scattering in fluorescence diffuse optical tomography}, Inverse Problems \&
  Imaging, 10 (2016), p.~227.

\bibitem{OLearyRoseberry_Villa_Chen_Ghattas_2022}
{\sc T.~O'Leary-Roseberry, U.~Villa, P.~Chen, and O.~Ghattas}, {\em
  Derivative-informed projected neural networks for high-dimensional parametric
  maps governed by {PDEs}}, Computer Methods in Applied Mechanics and
  Engineering, 388 (2022), p.~114199.

\bibitem{Plessix_2006}
{\sc R.-E. Plessix}, {\em A review of the adjoint-state method for computing
  the gradient of a functional with geophysical applications}, Geophysical
  Journal International, 167 (2006), pp.~495--503.

\bibitem{Pronzato_Walter_1988}
{\sc L.~Pronzato and E.~Walter}, {\em Robust experiment design via maximin
  optimization}, Mathematical Biosciences, 89 (1988), pp.~161--176.

\bibitem{Rainforth_Foster_Ivanova_BickfordSmith_2024}
{\sc T.~Rainforth, A.~Foster, D.~R. Ivanova, and F.~B. Smith}, {\em {Modern
  Bayesian Experimental Design}}, Statistical Science, 39 (2024), pp.~100 --
  114.

\bibitem{Rojas_Welsh_Goodwin_Feuer_2007}
{\sc C.~R. Rojas, J.~S. Welsh, G.~C. Goodwin, and A.~Feuer}, {\em Robust
  optimal experiment design for system identification}, Automatica, 43 (2007),
  pp.~993--1008.

\bibitem{Stuart_2010}
{\sc A.~M. Stuart}, {\em Inverse problems: A {Bayesian} perspective}, Acta
  Numerica, 19 (2010), p.~451–559.

\bibitem{Sunseri_Alexanderian_Hart_Waanders_2024}
{\sc I.~Sunseri, A.~Alexanderian, J.~Hart, and B.~v.~B. Waanders}, {\em
  Hyper-differential sensitivity analysis for nonlinear {B}ayesian inverse
  problems}, International Journal for Uncertainty Quantification, 14 (2024),
  p.~1–20.

\bibitem{Telen_Logist_VanDerlinden_VanImpe_2012}
{\sc D.~Telen, F.~Logist, E.~Van~Derlinden, and J.~F. Van~Impe}, {\em Robust
  optimal experiment design: A multi-objective approach}, IFAC Proceedings
  Volumes, 45 (2012), pp.~689--694.
\newblock 7th Vienna International Conference on Mathematical Modelling.

\bibitem{Telen_Vercammen_Logist_Van_Impe_2014}
{\sc D.~Telen, D.~Vercammen, F.~Logist, and J.~Van~Impe}, {\em Robustifying
  optimal experiment design for nonlinear, dynamic (bio)chemical systems},
  Computers \& Chemical Engineering, 71 (2014), pp.~415--425.

\bibitem{Ucinski_2005}
{\sc D.~Uci\'{n}ski}, {\em Optimal measurement methods for distributed
  parameter system identification}, Systems and Control Series, CRC Press, Boca
  Raton, FL, 2005.

\bibitem{Villa_Petra_Ghattas_2021}
{\sc U.~Villa, N.~Petra, and O.~Ghattas}, {\em {HIPPYlib: An Extensible
  Software Framework for Large-Scale Inverse Problems Governed by PDEs: Part I:
  Deterministic Inversion and Linearized Bayesian Inference}}, ACM Trans. Math.
  Softw., 47 (2021).

\bibitem{Scipy_2020}
{\sc P.~Virtanen, R.~Gommers, T.~E. Oliphant, M.~Haberland, T.~Reddy,
  D.~Cournapeau, E.~Burovski, P.~Peterson, W.~Weckesser, J.~Bright, S.~J. {van
  der Walt}, M.~Brett, J.~Wilson, K.~J. Millman, N.~Mayorov, A.~R.~J. Nelson,
  E.~Jones, R.~Kern, E.~Larson, C.~J. Carey, {\.I}.~Polat, Y.~Feng, E.~W.
  Moore, J.~{VanderPlas}, D.~Laxalde, J.~Perktold, R.~Cimrman, I.~Henriksen,
  E.~A. Quintero, C.~R. Harris, A.~M. Archibald, A.~H. Ribeiro, F.~Pedregosa,
  P.~{van Mulbregt}, and {SciPy 1.0 Contributors}}, {\em {{SciPy} 1.0:
  Fundamental Algorithms for Scientific Computing in Python}}, Nature Methods,
  17 (2020), pp.~261--272, \url{https://doi.org/10.1038/s41592-019-0686-2}.

\bibitem{Wald_1945}
{\sc A.~Wald}, {\em Statistical decision functions which minimize the maximum
  risk}, Annals of Mathematics, 46 (1945), p.~265–280.

\bibitem{Wu_Chen_Ghattas_2023}
{\sc K.~Wu, P.~Chen, and O.~Ghattas}, {\em A fast and scalable computational
  framework for large-scale high-dimensional bayesian optimal experimental
  design}, SIAM/ASA Journal on Uncertainty Quantification, 11 (2023),
  pp.~235--261.

\end{thebibliography}

\appendix

\section{Variational Tools}
\label{sec:variational-tools}
The gradient and the Hessian of the cost functional $\Phi$ \eqref{eq:postmean} are essential components of our proposed methods.
To compute these derivatives, we rely on adjoint-based gradient
computation~\cite{Plessix_2006}, derived using a formal Lagrangian approach.
Here, we outline the adjoint-based expressions for the gradient and Hessian
apply.
We begin by defining the Lagrangian
\begin{equation} \label{eq:lagrangian}
	\mc{L}(\state, \invparam, \test)
	= \frac{1}{2} \norm{\obs - \Obs\state}_{\noisecovinv}^2
	+ \frac{1}{2} \norm{\invparam - \priormean}_{\priorcovinv}^2
	+ \weakpde(\state, \invparam, \test) \,,
\end{equation}
where $\weakpde$ is the variational form of the forward problem \eqref{eq:weak-pde}.
In this context, $p$ is called the adjoint variable.
The gradient of $\Phi$ is given by
the variation of this Lagrangian with respect to $\invparam$, assuming the variations
of $\mc{L}$ with respect to $u$ and $p$ vanish. Namely,
\begin{subequations} \label{eq:gradient-system}
	\begin{equation} \label{eq:gradient}
		\mc{G}(\invparam)(\tilde{\invparam})
		= \inp{\tilde{\invparam}}{\invparam - \priormean}_{\priorcovinv}
		+ \inp{\tilde{\invparam}}{\weakpde_{\invparam}(\state, \invparam, \test)} \,,
	\end{equation}
	where $u$ and $p$ satisfy
	\begin{align}
		\label{eq:state}
		\inp{\tilde{\test}}{\weakpde_\test(\state, \invparam, \test)}
		 & = 0
		\quad \forall \tilde{\test} \in \testspace \,, \\
		\label{eq:adjoint}
		\inp{\tilde{\state}}{\weakpde_{\state}(\state, \invparam, \test)}
		+ \inp{\tilde{\state}}{\Obs^* \noisecovinv (\obs - \Obs\state)}
		 & = 0
		\quad \forall \tilde{\state} \in \testspace \,.
	\end{align}
\end{subequations}

To compute the Hessian action of $\Phi$, we follow a
Lagrange multiplier approach to differentiate through the gradient
\eqref{eq:gradient} constrained by the state and adjoint equations
\eqref{eq:state}, and \eqref{eq:adjoint}, respectively.
A detailed discussion of deriving the action of the Hessian via adjoint-based
techniques can be found in \cite{Villa_Petra_Ghattas_2021}.
As explained later, in our proposed approach we only require the data-misfit Hessian $\Hm$.
To derive the adjoint-based data-misfit Hessian action
we consider the meta-Lagrangian
\begin{multline*} 
	\mc{L}^H(\state, \invparam, \test, \hat{\state}, \hat{\invparam}, \hat{\test})
	=
	\\
	\inp{\hat{\invparam}}{\weakpde_{\invparam}(\state, \invparam, \test)}
	+ \inp{\hat{\test}}{\weakpde_\test(\state, \invparam, \test)}
	+ \inp{\hat{\state}}{\weakpde_{\state}(\state, \invparam, \test)}
	+ \inp{\hat{\state}}{\Obs^* \noisecovinv (\obs - \Obs\state)} \,.
\end{multline*}
Through a similar process as before, we take a variation of $\mc{L}^H$ with respect to
$\invparam$ and constrain it by letting variations of $\mc{L}^H$ with respect to
$\hat\state$ and $\hat \test$ vanish. This yields the expression for the data-misfit
Hessian action:
\begin{subequations}
	\begin{equation} \label{eq:data-misfit-hessian}
		\Hm(\invparam)(\hat{\invparam}, \tilde{\invparam})
		= \inp{\tilde{\invparam}}{
			\weakpde_{\invparam\invparam}(\state, \invparam, \test) \hat{\invparam}
			+ \weakpde_{\invparam\state}(\state, \invparam, \test) \hat{\state}
			+ \weakpde_{\invparam\test}(\state, \invparam, \test) \hat{\test}
		} \,,
	\end{equation}
	where for all $\tilde{\test} \in \testspace$ and $\tilde{\state} \in \testspace$,
	\begin{align}
		\label{eq:incr-state}
		\inp{\tilde{\test}}{\weakpde_{\test\state}(\state, \invparam, \test) \hat{\state}}
		+ \inp{\tilde{\test}}{\weakpde_{\test\invparam}(\state, \invparam, \test)
			\hat{\invparam}}
		 & = 0 \,, 
		\\
		\label{eq:incr-adjoint}
		\inp{\tilde{\state}}{
			\weakpde_{\state\test}(\state, \invparam, \test) \hat{\test}
			+ \weakpde_{\state\state}(\state, \invparam, \test) \hat{\state}
			+ \weakpde_{\state\invparam}(\state, \invparam, \test) \hat{\invparam}
		}
		+ \inp{\tilde{\state}}{\Obs^* \noisecovinv \Obs \hat{\state}}
		 & = 0 \,. 
	\end{align}
\end{subequations}
Finally, the Gauss-Newton Hessian is obtained by dropping the terms involving
the adjoint variable; see~\cite{Bangerth_2008} for details.
In the present setting, the Gauss-Newton data-misfit Hessian action is given by
\begin{subequations}
	\begin{equation}
		\mc{H}(\invparam)(\hat{\invparam}, \tilde{\invparam})
		= \inp{\tilde{\invparam}}{
			\weakpde_{\invparam\invparam}(\state, \invparam, \test)
			+ \weakpde_{\invparam\state}(\state, \invparam, \test) \hat{\state}
			+ \weakpde_{\invparam\test}(\state, \invparam, \test) \hat{\test}
		},
	\end{equation}
	where
	for all $\tilde{\test} \in \testspace$
	and $\tilde{\state} \in \testspace$,
	\begin{align}
		\inp{\tilde{\test}}{\weakpde_{\test\state}(\state, \invparam, \test) \hat{\state}}
		+ \inp{\tilde{\test}}{\weakpde_{\test\invparam}(\state, \invparam, \test)
			\hat{\invparam}}
		 & = 0 \,, 
		\\
		\inp{\tilde{\state}}{
			\weakpde_{\state\test}(\state, \invparam, \test) \hat{\test}
		}
		+ \inp{\tilde{\state}}{\Obs^* \noisecovinv \Obs \hat{\state}}
		 & = 0 \,. 
	\end{align}
\end{subequations}

\iftoggle{arxiv}{
	\null \vfill
	\begin{flushright}
		\scriptsize \framebox{\parbox{5.5in}
			{
				The submitted manuscript has been created by UChicago Argonne, LLC, Operator of
				Argonne National Laboratory (``Argonne"). Argonne, a U.S. Department of Energy
				Office of Science laboratory, is operated under Contract No. DE-AC02-06CH11357. The
				U.S. Government retains for itself, and others acting on its behalf, a paid-up
				nonexclusive, irrevocable worldwide license in said article to reproduce, prepare
				derivative works, distribute copies to the public, and perform publicly and display
				publicly, by or on behalf of the Government. The Department of Energy will provide
				public access to these results of federally sponsored research in accordance with
				the DOE Public Access Plan. http://energy.gov/downloads/doe-public-access-plan.
			}
		}
		\normalsize
	\end{flushright}
}{}

\end{document}